\documentclass[times]{elsarticle}
\usepackage{tabularx}
\usepackage{color}
\usepackage{moreverb} 
\usepackage[colorlinks,bookmarksopen,bookmarksnumbered,citecolor=red,urlcolor=red,pdfauthor=author]{hyperref}
\usepackage{amsfonts,amsmath,amssymb}
\usepackage{psfrag,pstricks,pst-node}

\usepackage[section]{placeins}
\usepackage{float}
\usepackage{subfigure} 
\usepackage{mathrsfs}
\usepackage{graphicx}
\usepackage[ruled,vlined]{algorithm2e}
\usepackage{fancyhdr}  
\usepackage{psfrag} 
\usepackage[export]{adjustbox} 
\usepackage{enumerate} 
\usepackage{comment}
\usepackage{lscape}
\usepackage{appendix}
\usepackage{algorithmic}
\usepackage{bm}
\usepackage{booktabs}
\usepackage{threeparttable}
\graphicspath{{./figs/}{./model Figure/}}
\DeclareGraphicsExtensions{.pdf,.eps,.PNG,.png}

\newcommand\BibTeX{{\rmfamily B\kern-.05em \textsc{i\kern-.025em b}\kern-.08emT\kern-.1667em\lower.7ex\hbox{E}\kern-.125emX}}

\usepackage{pgfplots}
\pgfplotsset{compat=1.13}

\usepackage{geometry}
\usepackage{graphicx,array}
\usepackage{subcaption} 
\usepackage{amsmath}
\usepackage[numbers]{natbib}
\usepackage{url}
\usepackage{amsmath}
\usepackage{threeparttable} 
\usepackage{multirow}
\usepackage{changepage}
\usepackage{booktabs}
\usepackage{fullpage}
\usepackage{times}
\usepackage{fancyhdr,graphicx,amsmath,amssymb}
\usepackage[ruled,vlined]{algorithm2e}
\usepackage{float}
\usepackage[T1]{fontenc}
\usepackage[utf8]{inputenc}
\usepackage{bm}
\usepackage{booktabs}
\usepackage{nomencl}
\usepackage{enumitem}
\usepackage{makecell}
\usepackage{arydshln}

\usepackage{makecell} 

\begin{document}
 
\begin{frontmatter}
\renewcommand{\thefootnote}{\fnsymbol{footnotemark}}
 
\fancypagestyle{plain}{%
\fancyhf{} 
\fancyhead[RO,RE]{\thepage} 
}

\title{Nonlinear Model Reduction by Probabilistic Manifold Decomposition}   
      \author[lab1]{Jiaming Guo}
    \author[lab1]{Dunhui Xiao\corref{cor1}}     
         \cortext[cor1]{Corresponding author}
    \ead{xiaodunhui@tongji.edu.cn}

    \address[lab1]{School of Mathematical Sciences,
Key Laboratory of Intelligent Computing and Applications(Ministry of Education), Tongji University, Shanghai 200092, CHINA}

\begin{abstract}
This paper presents a novel non-linear model reduction method: Probabilistic Manifold Decomposition (PMD), which provides a powerful framework for constructing non-intrusive reduced-order models (ROMs) by embedding a high-dimensional system into a low-dimensional probabilistic manifold and predicting the dynamics. Through explicit mappings, PMD captures both linearity and non-linearity of the system. A key strength of PMD lies in its predictive capabilities, allowing it to generate stable dynamic states based on embedded representations.

The method also offers a mathematically rigorous approach to analyze the convergence of linear feature matrices and low-dimensional probabilistic manifolds, ensuring that sample-based approximations converge to the true data distributions as sample sizes increase. These properties, combined with its computational efficiency, make PMD a versatile tool for applications requiring high accuracy and scalability, such as fluid dynamics simulations and other engineering problems. By preserving the geometric and probabilistic structures of the high-dimensional system, PMD achieves a balance between computational speed, accuracy, and predictive capabilities, positioning itself as a robust alternative to the traditional model reduction method.

\end{abstract}

\begin{keyword} 
 PMD, Model order reduction; Manifold learning; Geodesic distance; Non-intrusive ROM
\end{keyword}
\end{frontmatter}

\section{Introduction}
\label{sec:intro}
\vspace{-2pt}
In computational engineering, solving high-dimensional problems has been a critical challenge due to the large dimension size resulting from the discretization of the governing equations of complex systems \cite{diez2021nonlinear}. These high-dimensional spaces, where the degrees of freedom correspond to the discretization of the equations, are computationally expensive. However, it is frequently observed that the solution sets of these problems lie on low-dimensional manifolds, governed by a smaller set of independent parameters. This observation highlights an opportunity for reduced-order modeling (ROM), a technique designed to alleviate the computational burden by reducing both dimensionality and degrees of freedom. ROMs offer a pathway to more efficient solutions while preserving the essential features of the original system. By projecting high-dimensional data onto lower-dimensional manifolds, ROMs not only reduce computational complexity but also retain key dynamics, enabling faster simulations and real-time solutions. This makes ROM indispensable in time-critical applications, such as parameter optimization, sensitivity analysis, and large-scale simulations.

The use of reduced-order models has become indispensable across various domains, enabling efficient solutions without sacrificing the fidelity of the model. By leveraging lower-dimensional representations, ROMs make it possible to model complex systems efficiently, without losing the accuracy necessary for high-fidelity simulations. Recent advances in machine learning have further enhanced the capabilities of ROMs, particularly in improving their predictive accuracy and generalization. The integration of machine learning techniques, such as neural networks, has significantly improved predictive modeling in fields such as computational fluid dynamics (CFD)\cite{xie2018data}, climate prediction \cite{o2018using}, engineering design \cite{willard2022integrating}, and biomedical engineering \cite{foster2014machine}. These developments allow for more accurate models and faster simulations, tackling the complexities of real-world problems.


Existing model reduction methods include balanced truncation\cite{heinkenschloss2008balanced}, the reduced-basis method\cite{peterson1989reduced}, rational interpolation\cite{baur2011interpolatory}, and the proper orthogonal decomposition (POD)\cite{berkooz1993proper}. Within these methods, the ROM restricts the state to evolve in a linear subspace, which imposes a fundamental limitation to the accuracy of the resulting ROM\cite{lee2020model}. While POD has been widely applied across various fields due to its effectiveness in reducing dimensionality, it struggles to handle the nonlinearities inherent in many real-world systems\cite{schlegel2015long,chaturantabut2010nonlinear,osth2014need}. Additionally, methods such as the empirical interpolation method (EIM)\cite{barrault2004empirical}, the discrete empirical interpolation method (DEIM)\cite{peherstorfer2014localized}, a combination of quadratic expansion and DEIM (residual DEIM)\cite{xiao2014non}, and Gauss-Newton methods with approximated tensors (GNAT)\cite{carlberg2013gnat} have been developed to mitigate the nonlinear inefficiencies associated with POD. While the combination of DEIM with POD effectively improves computational efficiency and stability in reducing nonlinear terms, it does not address the fundamental limitation of POD in capturing the accuracy of nonlinear dynamics.

To address these limitations, various nonlinear dimensionality reduction methods have been proposed \cite{zhang2004principal}, including isometric mapping (Isomap) \cite{balasubramanian2002isomap}, autoencoders \cite{fu2023non,lee2020model}, local linear embedding (LLE) \cite{roweis2000nonlinear}, t-SNE \cite{van2008visualizing}, and diffusion maps \cite{coifman2006diffusion}. By doing so, they offer a more accurate representation of the system's underlying dynamics and are particularly effective in revealing the complex nonlinear relationships that arise in real-world systems. However, despite their advantages, these methods also have limitations. One of the primary challenges is their poor interpretability, as they often transform the data into a lower-dimensional space that is not easily understood or linked to the physical or governing variables of the system. Furthermore, many of these techniques struggle to capture the global structure of the data adequately, focusing more on local relationships and neglecting the broader, more holistic dynamics of the system.

In this paper, we propose a novel non-linear model reduction method, Probabilistic Manifold Decomposition (PMD). By leveraging existing solvers, PMD efficiently constructs nonlinear manifolds without modifying the high-fidelity model \cite{xiao2015non}, making it non-intrusive and applicable to a wide range of scientific and engineering problems. 

This method itself is relatively new and combines manifold learning with probabilistic modeling to form a low-dimensional representation that accurately approximates the full system. However, due to the high nonlinearity of the Navier-Stokes equations, It is difficult for us to preserve the main features of the system. To mitigate this problem, PMD method uses two distinct steps: first, it reduces the system’s dimensionality through the linear component extraction via SVD, followed by projecting the nonlinear components onto a probabilistic manifold. The main linear features are captured using Singular Value Decomposition (SVD), while the main nonlinear features are embedded onto a low-dimensional manifold using a probabilistic approach.

Once the linear terms are captured through a reduced basis, the nonlinear terms are treated separately. The method identifies the nonlinear residual by subtracting the linear approximation from the original data. These residuals are then mapped onto a probabilistic manifold that preserves both local and global geometric relationships in the system’s dynamics. This step allows PMD to capture the inherent nonlinearities in the system while keeping the dimensionality manageable.

To construct this manifold, PMD uses a Markov process, which helps to create a similarity graph between the data points based on how close they are in the system’s state space. This graph is built by calculating the similarity between data points, capturing both local and global relationships in the data.

Once the similarity graph is constructed, PMD uses a process that simulates a random walk across the data points, which generates a transition matrix. This matrix represents the probabilities of transitioning from one state to another in the system, based on the graph of similarities. By applying this transition matrix over multiple steps, PMD captures the system's dynamics and uncovers the low-dimensional manifold that accurately represents the system's behavior.

A key feature of the PMD method is the use of geodesic distance to measure the intrinsic geometry of the manifold. This ensures that the nonlinear components are correctly represented, unlike traditional methods that rely solely on Euclidean distance. The geodesic distance is computed using the Floyd-Warshall algorithm on the weighted adjacency graph, which represents the similarities between the data points. This process is essential for preserving the manifold's structure, especially in highly nonlinear systems like the Navier-Stokes equations.

Once the manifold is constructed, PMD predicts the future states of the system using a two-step process. For linear components, the method draws on the principles of dynamic mode decomposition (DMD) in the reduced space to predict their evolution. For nonlinear components, PMD evolves the probabilistic manifold over time, maintaining the low-dimensional representation of the system's nonlinear dynamics. The two predicted components are then mapped back to the original high-dimensional space using a lift mapping, which reconstructs the full state of the system. This ensures that both the linear and nonlinear dynamics are captured accurately, providing a robust and computationally efficient reduced-order model.

The convergence of PMD is another significant strength. The method ensures that the reduced model converges to the true dynamics of the full system as more data becomes available. Through careful analysis, it is demonstrated that both the linear and nonlinear features converge at an appropriate rate, meaning that the reduced model can approximate the full model’s behavior with increasing accuracy as the data set grows. This convergence is crucial for ensuring that PMD remains effective in real-world applications, where the system’s dynamics may evolve over time.

The capabilities of PMD are demonstrated via two fluid problems: flow past a cylinder with Reynolds numbers $Re=300$ and $Re=5000$, and the lock exchange problem. The results show the PMD accurately captures vortex shedding and diffusion processes.

The structure of this paper is as follows. Section \ref{sec:govern} introduces the theoretical background and governing equations, including the Navier-Stokes equations. Section \ref{sec:method} details the PMD framework, focusing on its approach to dimensionality reduction, prediction, and integration with data-driven regression methods. Section \ref{sec:convergence} presents the convergence analysis of PMD, establishing its theoretical robustness and accuracy. Section \ref{sec:experiments} provides numerical examples that demonstrate PMD's effectiveness in fluid dynamics, including case studies on flow around a cylinder and the lock exchange problem. Finally, Section \ref{sec:summary} offers conclusions and discusses future directions for advancing PMD-based reduced-order modeling in complex systems.

\section{Governing equations}
\label{sec:govern}
The equations are non-compressible Navier–Stokes equations describing the conservation of mass and momentum for fluids:
\begin{equation}\label{eq:equation1}
\nabla \cdot \bf{u} = 0.
\end{equation}
Equation \ref{eq:equation1} is the mass conservation equation for fluids, commonly known as the continuity equation, $\nabla \cdot \bf{u}$ represents the divergence of the velocity field $\bf{u}$.
\begin{equation}\label{eq:equation2}
\frac{\partial \mathbf{u}}{\partial t} + (\mathbf{u} \cdot \nabla) \mathbf{u} + f \mathbf{k} \times \mathbf{u} = -\nabla p + \nabla \cdot \tau.
\end{equation}
Equation \ref{eq:equation2} is the momentum conservation equation for fluids. The vector $\bf{u}$ represents the velocity. $p: = p/\rho_0$ refers to the modified pressure, where $\rho_0$ is the constant reference density of the fluid.
$f\bf{k} \times \bf{u}$ describes the effect of the Coriolis force on fluid motion, where $\it{f}$ is the Coriolis parameter, usually related to the rotational angular velocity of the Earth, and $\bf{k}$ is the unit vector along the rotation axis of the Earth. $\tau$ represents the stress tensor, describing the internal stress distribution within the fluid. The stress tensor is usually expressed as:
\begin{equation}
\tau := \mu (\nabla\mathbf{u} + (\nabla\mathbf{u})^T).
\end{equation}

\section{Probabilistic manifold decomposition (PMD)}
\label{sec:method}

PMD is a data-driven nonlinear model reduction method that efficiently maps a high-dimensional system to a low-dimensional manifold while preserving its essential structure.  PMD eliminates redundancies and captures both linear and nonlinear features from the system, making it effective for complex dynamic systems where traditional linear methods may fail. This approach combines probabilistic modelling with dimensionality reduction techniques, enabling an accurate approximation of the intrinsic geometry of the dynamic states and facilitating efficient prediction.

\subsection{Linear dimensionality reduction in PMD}
The linear dimensionality reduction method is achieved via the traditional SVD method as it proved to be an efficient method and is widely used in engineering.  
\begin{equation}\label{eq:svd}
U = 
\begin{bmatrix}
\mathbf{u}_1 & \mathbf{u}_2 & \cdots & \mathbf{u}_m \\
\end{bmatrix},\quad
\tilde{U} = 
\begin{bmatrix}
\frac{\mathbf{u_1} - \mu_1}{\overline{\sigma}_1} & \cdots & \frac{\mathbf{u_m} - \mu_m}{\overline{\sigma}_m} 
\end{bmatrix},
\quad
\tilde{U} = Q \Sigma V^.
\end{equation}
Here, \( U \in R^{n \times m} \) is the snapshot of simulation (or experimental) data, where \( n \) denotes the number of spatial points, and \( m \) is the number of snapshots. \( \mu_i \) and \( \overline{\sigma}_i \) represent the mean and standard deviation of \(\mathbf{u}_i\) respectively. \(Q \in \mathbb{R}^{n \times n}\) and \(V \in \mathbb{R}^{m \times m}\) are orthogonal matrices, and \(\Sigma \in \mathbb{R}^{n \times m}\) contains singular values \(\sigma_1 \geq \sigma_2 \geq \dots \geq \sigma_{\min(n, m)}\).

To retain the primary variance, PMD selects the first \(r\) singular values, ensuring:
\begin{equation}
\sum_{i=1}^r \sigma_i > (1-\epsilon) \sum_{i=1}^m \sigma_i,
\end{equation}
The reduced representation is given by:
\begin{equation}\label{eq:linear}
\tilde{U}_r = Q_r \Sigma_r V_r^T,
\end{equation}
where \(Q_r\in\mathbb{R}^{n \times r}\), \(\Sigma_r\in\mathbb{R}^{r \times r}\), and \(V_r\in\mathbb{R}^{m \times r}\) are truncated matrices containing the first \(r\) components. The reconstruction error is:
\begin{equation}
||X||_F^2 = \sum_{i=r+1}^m \sigma_i^2,
\end{equation} 
\begin{equation}\label{eq:error}
X := \tilde{U} - P\tilde{U} \in R^{n\times m},
\end{equation}
represents the components of the normalized matrix $\tilde{U}$ that are not captured by the linear basis. The matrix $P$ is an orthogonal projection operator, defined as:
\begin{equation}
P := V_r  (V_r^TV_r) ^{-1}V_r^T.
\end{equation}
Since $V_r$ is an orthonormal matrix, we have $V_r^TV_r = I_r \in R^{r \times r}$. Thus, equation \ref{eq:error} can be reformulated as:
\begin{equation}\label{eq:error1}
X := (I - V_r V_r^T) \tilde{U} \in R^{n\times m}.
\end{equation}


\subsection{Nonlinear probabilistic manifold construction in PMD}\label{sec:3.2}

PMD constructs a data-driven nonlinear manifold that captures features beyond the linear basis. Inspired by the ideas of Taylor expansion, quadratic manifold \cite{geelen2023operator}, error analysis \cite{xiao2019error}, and diffusion maps \cite{coifman2006diffusion}, PMD uses a Markov process \cite{ethier2009markov} with transition probabilities to build a low-dimensional probabilistic manifold embedding, preserving both local and global geometric relationships in the dynamics. 

\subsubsection{Construction of the weighted adjacency graph}\label{sec:3.2.1}

As shown in equation \ref{eq:error1}, after obtaining the projection error (or residual) $X$ of the nonlinear components of the original data matrix through error analysis, PMD aims to learn a low-dimensional probabilistic model of this residual matrix. To achieve this, PMD first constructs a weighted adjacency graph $W$ to represent the similarity between data points, where the similarity metric must be symmetric and positive definite. Specifically, this means: $W(i, j) = W(j, i), \quad  W(i, j) \geq 0$. The similarity metric should be related to the distance between data points, specifically:
\begin{equation}
W(i, j) := K(||\mathbf{x}_i - \mathbf{x}_j||,\epsilon), \quad i, j = 1, 2, \cdots, m.
\end{equation}
Here, $W \in R^{m \times m}$ is the weights,  $\mathbf{x}_i\in X$ and $K$ represents the kernel function selected for weight calculation. The parameter $\epsilon$ is related to the intrinsic properties of the dynamic system, while $||\cdot||$ denotes the chosen metric, which is closely related to the underlying geometric structure of $X$ for dimensionality reduction. The Gaussian kernel function is used for calculating the weight matrix W. 

\begin{equation} \label{eq:Gauss}
W  (i, j)  = exp  (-\frac{||\mathbf{x}_i - \mathbf{x}_j||^2}{\epsilon^2}) , \quad i, j = 1, 2, \cdots, m.
\end{equation}

The parameter $\epsilon$ controls the weight calculation by defining the region where the similarity metric is effective. It is related to the data's geometric structure: smaller $\epsilon$ values are preferred for complex, nonlinear, low-dimensional data, while larger $\epsilon$ values are better suited for sparse datasets, ensuring appropriate neighborhood size.

\subsubsection{Geodesic distance learning of the probabilistic manifold}
\begin{figure}[tbhp]
    \centering
    \includegraphics[height=2cm]{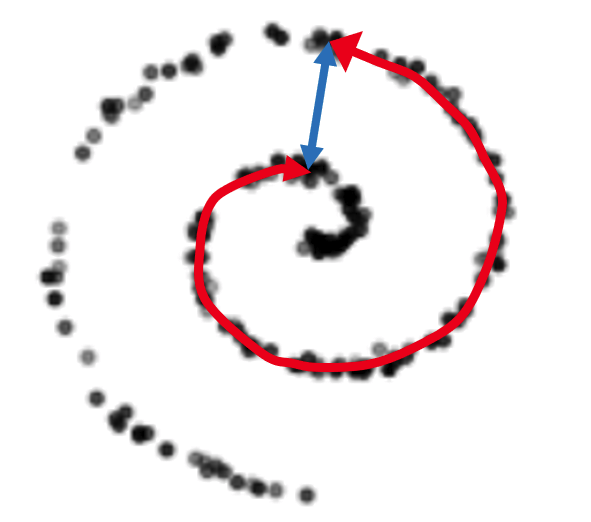}
    \caption{Euclidean distance (blue) vs. geodesic distance (red) along the manifold. The Euclidean distance fails to capture the manifold's intrinsic geometry, whereas the geodesic distance reflects the true relationships.}
    \label{figure:1}
\end{figure}

In traditional nonlinear dimensionality reduction methods, such as diffusion maps and Locally Linear Embedding (LLE), the Euclidean distance is commonly used to compute pairwise distances and determine weights between points, as shown in \ref{sec:3.2.1}. However, this approach fails to accurately capture the intrinsic geometry of the data manifold. In this work, the PMD method uses geodesic distance instead. This ensures a more faithful representation of the "true distance" along the manifold, providing a more accurate depiction of the data's structure, as illustrated in \ref{figure:1}.

As illustrated in \ref{figure:1}, the blue arrow represents the Euclidean distance, which measures the straight-line distance in the embedding space. However, this can introduce significant errors, especially for data lying on a curved manifold. The red arrow represents the geodesic distance, which accounts for the curvature and intrinsic geometry of the manifold, leading to a more faithful representation of the data relationships.

To compute the geodesic distance \cite{bronstein2006efficient}, PMD employs the Floyd-Warshall algorithm on a weighted adjacency graph $W$. This graph is constructed such that the edges represent pairwise distances between neighboring points. The geodesic distance is then approximated by finding the shortest paths between all pairs of points, as described by the following formula:
\begin{equation}
    d_g(\mathbf{x}_i, \mathbf{x}_j) = \min_{k} \left( d_g(\mathbf{x}_i, \mathbf{x}_k) + d_g(\mathbf{x}_k, \mathbf{x}_j) \right),
\end{equation}
where $d_g(\mathbf{x}_i, \mathbf{x}_j)$ denotes the shortest path distance between points $\mathbf{x}_i$ and $\mathbf{x}_j$, and $k$ iterates over all intermediate vertices.
The Floyd-Warshall algorithm iteratively updates the shortest path estimates using the recursive rule:
\begin{equation}
    d_g^{(k)}(\mathbf{x}_i, \mathbf{x}_j) = \min \left( d_g^{(k-1)}(\mathbf{x}_i, \mathbf{x}_j), d_g^{(k-1)}(\mathbf{x}_i, \mathbf{x}_k) + d_g^{(k-1)}(\mathbf{x}_k, \mathbf{x}_j) \right),
\end{equation}
where $d_g^{(k)}(\mathbf{x}_i, \mathbf{x}_j)$ represents the shortest path distance considering up to the $k$-th vertex as an intermediate point. Initially, the distances $d_g^{(0)}(\mathbf{x}_i, \mathbf{x}_j)$ are set as the edge weights, $W(i,j)$, for direct neighbors. Once the algorithm converges, the resulting matrix $D = [d_g(\mathbf{x}_i, \mathbf{x}_j)]$ contains the geodesic distances between all pairs of points.  

\subsubsection{Construction of the transition matrix}

\begin{figure}[tbhp]
 \centering
 \includegraphics[height=4cm]{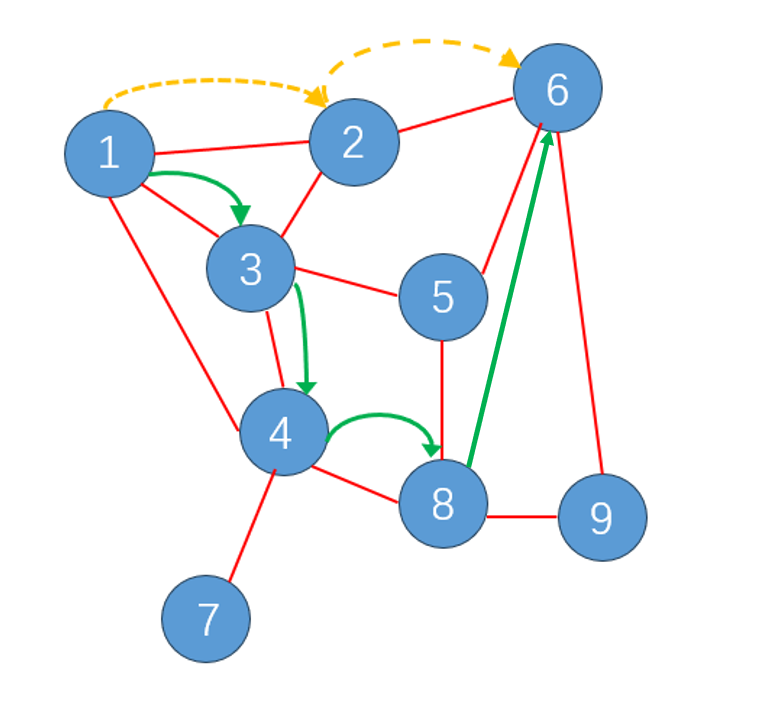}
 \caption{Random walk process}\label{fig:2}
\end{figure}

After calculating the weights, PMD normalizes the adjacency matrix to construct a Markov transition matrix $P \in R^{m \times m}$. This matrix represents the transition probabilities between data points, where closer points have higher probabilities, and distant points have lower probabilities, based on the chosen distance metric.
\begin{equation} \label{eq:sample1}
P(i, j) = \frac{W(i, j)}{\sum_{k=1}^m W(i, k)}, \quad i, j = 1, 2, \cdots, m.
\end{equation}
The matrix \( P \) obtained above is the one-step transition matrix in a Markov process. To calculate the probability of transition along a path, as shown in \ref{fig:2}, we need to compute transition matrices in the following steps. For example, the two-step transition matrix is used to determine the transition probability along the yellow arrow, while the four-step transition matrix is required for the green arrows.
The \( t \)-step transition matrix is given by:
\begin{equation}
    \tilde{P} = P^t.
\end{equation}

In the PMD method, selecting the optimal number of transition steps $t$ is crucial to capture the local and global structures of the data manifold. More steps emphasize global patterns by smoothing small variations, while fewer steps preserve local details and subtle variations.

Finding the right balance in $t$ ensures that the low-dimensional representation accurately reflects the data manifold structure, capturing both fine local features and global geometric relationships. This balance enhances the accuracy and efficiency of the embedding.

\subsubsection{Construction of the probabilistic manifold}

Once the transition step \( t \) is determined, PMD performs eigenvalue decomposition of the Markov transition matrix \( P \), extracting the leading eigenvalues and eigenvectors to form a low-dimensional embedding. This embedding preserves the geometric structure of the data while avoiding the computational cost of high-dimensional distance calculations. 

PMD identifies minimal coordinates for efficient representation via the mapping \( \Phi: \mathbb{R}^n \to \mathbb{R}^r \) from the high dimensional space into low dimensional manifold space:
\begin{equation} \label{eq:embedding}
\Phi (\mathbf{x}_i) = 
\begin{bmatrix}
\lambda_1^t \cdot \phi_1 (x_i) \\
\lambda_2^t \cdot \phi_2 (x_i) \\
\vdots \\
\lambda_r^t \cdot \phi_r (x_i) \\
\end{bmatrix}, \qquad i = 1, 2, \cdots, r,
\end{equation}
where \( \lambda_0 \ge \lambda_1 \ge \cdots \ge \lambda_r \ge 0 \) are eigenvalues of \( P \), and \( \phi_j \) are the corresponding eigenvectors. The component \( \phi_j(x_i) \) denotes the \( i \)-th component of the eigenvector associated with \( \lambda_j \). Since \( P \) is a transition matrix, \( \lambda_0 = 1 \) and \( \phi_0 = [1, 1, \cdots, 1]^T \).

The transition step \( t \) controls eigenvalue behavior, with \( \lambda_i^t \) decreasing as \( t \) increases. A gap between \( \lambda_r \) and \( \lambda_{r+1} \) ensures \( \lambda_{r+1}^t \to 0 \) while \( \lambda_r^t \) remains nonzero, retaining essential geometric information. 

This yields a low-dimensional Euclidean representation \( \Phi (\mathbf{x}_i) \) approximating true distances in the original space. PMD reduces dimensionality from \( \mathbb{R}^{n \times m} \) to \( \mathbb{R}^{r \times m} \), producing the embedding matrix:
\begin{equation}\label{eq:manifold}
\Phi = 
\begin{bmatrix}
\Phi (\mathbf{x}_1) & \Phi (\mathbf{x}_2) & \cdots & \Phi (\mathbf{x}_m) \\
\end{bmatrix}.
\end{equation}

\subsubsection{Construction of the lift mapping}
The lift mapping involves projecting the low-dimensional manifold back into the high-dimensional space. It constructs a mapping \(\psi: \mathbb{R}^r \to \mathbb{R}^n\) that reconstructs a high-dimensional vector \(\mathbf{x}\) from its low-dimensional representation \(\hat{\mathbf{x}}\):  
\begin{equation}\label{eq:lift}
\psi(\hat{\mathbf{x}}) = \mathbf{x}.
\end{equation}

To compute \(\psi\), we utilize \(K\)-nearest neighbors from the dataset, \(\{\mathbf{x}_j\}_{j=1}^K \in \mathbb{R}^n\) and their low-dimensional counterparts \(\{\hat{\mathbf{x}}_j\}_{j=1}^K \in \mathbb{R}^r\), and perform polynomial interpolation on the probabilistic manifold. The matrices are expressed as:  
\begin{equation}
\hat{X} = \begin{bmatrix}
 \mathbf{\hat{x}}_1 & \mathbf{\hat{x}}_2 & \cdots & \mathbf{\hat{x}}_K 
\end{bmatrix} \in \mathbb{R}^{r \times K}, \quad X = \begin{bmatrix}
 \mathbf{x}_1 & \mathbf{x}_2 & \cdots & \mathbf{x}_K 
\end{bmatrix} \in \mathbb{R}^{n \times K}.
\end{equation}
The interpolation coefficients are obtained by minimizing:  
\begin{equation}
K = \underset{K \in \mathbb{R}^{n \times r}}{\arg\min} \|K\hat{X} - X\|_F^2,
\end{equation}
where regularization is applied to stabilize the solution.  
To address ill-conditioning and overfitting, we adopt regularization techniques such as L2 regularization (ridge regression) \cite{cortes2012l2} and kernel ridge regression \cite{vovk2013kernel}, which are particularly effective for capturing non-linear dependencies:  
\begin{equation}\label{eq:Ridge}
K = \underset{K\in \mathbb{R}^{n \times r}}{\arg\min} \left(\|K\hat{X} - X\|_F^2 + \lambda \|K\|_F^2\right),
\end{equation}
where \(\lambda\) is the regularization parameter.
For non-linear data, we utilize kernel ridge regression with a polynomial kernel to approximate the high-dimensional features.
\begin{equation}\label{eq:kernal}
K = \underset{K \in \mathbb{R}^{K \times K}}{\arg\min} \left( \| K\hat{X}_1 - X_1 \|_F^2 + \lambda \| K \|_F^2 \right).
\end{equation}
The kernel function is defined as:  
\begin{equation}\label{eq:poly}
\hat{X}_1 = (\hat{X}^T \hat{X} + c)^d, \quad X_1 = (X^T X + c)^d,
\end{equation}
where \(d\) is the degree of the polynomial kernel, and \(c\) is a tuning parameter. This allows PMD to accurately capture non-linear relationships and reconstruct the high-dimensional data effectively. It also balances fitting accuracy and model complexity, ensuring that critical features of the high-dimensional space are preserved while avoiding overfitting.

\subsection{Linear term prediction}
PMD divides the matrix $\overline{U}_{r} = \Sigma_rV_r^T\in \mathbb{R}^{r\times m}$ derived from equation \ref{eq:linear} into sequential matrix $\overline{U}_{r_1}$ and time-shifted matrix $\overline{U}_{r_2}$:
\begin{equation}\label{eq: devide}
   \overline{U}_{r_1} = [\overline{\mathbf{u}}(1), \overline{\mathbf{u}}(2), \ldots, \overline{\mathbf{u}}(m-1)], \quad 
   \overline{U}_{r_2} = [\overline{\mathbf{u}}(2), \overline{\mathbf{u}}(3), \ldots, \overline{\mathbf{u}}(m)].
\end{equation}

Then, we seek to explore the time-series relationship between the sequential matrix $\overline{U}_{r_1}$ and the time-shifted matrix $\overline{U}_{r_2}$ by constructing a regularized dynamic projection operator:
\begin{equation}
    \overline{U}_{r_2} = A_1\overline{U}_{r_1},
\end{equation}
the linear relationship between these matrices is captured via:
\begin{equation}
    A_1 = \underset{A_1 \in \mathbb{R}^{r \times r}}{\arg\min} \left(||A_1\overline{U}_{r_1} - \overline{U}_{r_2}||_F^2 + \lambda||A_1||_F^2\right),
\end{equation}
\begin{equation}
A_1 =  \overline{U}_{r_2} \overline{U}_{r_1}^T\left( \overline{U}_{r_1} \overline{U}_{r_1}^T+\lambda I\right)^{-1}
\end{equation}
where \( A_1 \) estimates the transformation, $\lambda$ is the regularization parameter, which prevents ill-conditioned matrix inversion. Then we perform spectral decomposition on the matrix $A_1$ to extract the PMD mode,
\begin{equation}
A_1 Z = DZ,
\end{equation}
where $Z$ is the eigenvector matrix, $D$ is the eigenvalue matrix, $z_j, d_j$ are the eigenvectors and eigenvalues. Given the initial condition \(\overline{\mathbf{u}}(i)\) on the low-dimensional linear manifold, we compute its projection in the PMD mode space:
\begin{equation}\label{eq:mode}
B := (Z^TZ + \lambda I)^{-1}Z^T\overline{\mathbf{u}}(i),
\end{equation}
where $\lambda$ is the regularization parameter shown above to improve numerical stability and \(B \in \mathbb{R}^r\).

The next time-level dynamic states are predicted via:
\begin{equation}
 \hat{\mathbf{u}}(i + k) = \sum_{j=1}^n b_j \lambda_j^k z_j,
\end{equation}
where \(b_j\) are the components of \(B\), and \(k\) denotes the future time step.

\subsection{Nonlinear term prediction}

After predicting linear features, PMD captures nonlinear dynamics by refining the low-dimensional representation obtained earlier. 
Inspired by \cite{evangelou2023double}, this models complex non-linear dynamics using the probabilistic manifold $\overline{\Phi} \in \mathbb{R}^{r \times m}$. The dynamics are expressed as:
\begin{equation}\label{eq:g}
\phi_{m+1} = g(\phi_1, \phi_2, \dots, \phi_m),
\end{equation}
where \( g \) represents the evolution function in the manifold.  The steps of obtaining the function \( g \) are:

\textbf{Step 1: Kernel construction}.
Constructing a Gaussian kernel to capture similarities between manifold points:
\begin{equation}
K_{\Phi}(\phi_i, \phi_j) = \exp\left(-\frac{\|\phi_i - \phi_j\|^2}{\varepsilon^2}\right),
\end{equation}
where $\varepsilon$ is the bandwidth parameter, $\|\phi_i - \phi_j\|$ denotes the 'true distance' we denote in \ref{sec:3.2} between points $\phi_i$ and $\phi_j$ in the low-dimensional probabilistic manifold space.

\textbf{Step 2: PMD operator and eigenvalue decomposition.}  
Normalizing the kernel and computing eigenvectors \( v_j \):
\begin{equation}
P_{\Phi} = D_{\Phi}^{-1} K_{\Phi}, \quad P_{\Phi} v_j = \sigma_j v_j,
\end{equation}
where $D_{\Phi}$ is the diagonal degree matrix with elements $D_{\Phi}(i,i) = \sum_{j=1}^m K_{\Phi}(\phi_i, \phi_j)$. The \( \sigma_j \) are the eigenvalues, satisfying \( 1 = \sigma_1 \geq \sigma_2 \geq \cdots \geq \sigma_r \geq 0 \), and \( v_j \) are the corresponding eigenvectors.

\textbf{Step 3: Construction of the target function.} 
In this step, the data points \( \phi_1, \phi_2, \ldots, \phi_m \) are divided into two subsets:
\begin{equation}
\Phi_1 = [\Phi(x_1), \ldots, \Phi(x_{m-1})], \quad \Phi_2 = [\Phi(x_2), \ldots, \Phi(x_m)].
\end{equation}
Then the ridge regression is used to fit a mapping \( \overline{W} \) such that:
\begin{equation}
\Phi_2 = \overline{W} \Phi_1,
\end{equation}
by minimizing:
\begin{equation}
\mathcal{L}(\overline{W}) = \| \Phi_2 - \overline{W} \Phi_1 \|_F^2 + \lambda \| \overline{W} \|_F^2.
\end{equation}
The closed-form solution is:
\begin{equation}
\overline{W} = \Phi_2 \Phi_1^\top (\Phi_1 \Phi_1^\top + \lambda I)^{-1}.
\end{equation}
After that, the next time step is predicted as:
\begin{equation}
\Phi(x_{m+1}) = \overline{W} \Phi(x_m).
\end{equation}

\textbf{Step 4: Projection of the target function.}
Projecting \( g(\Phi) \) (represented by \( \overline{W} \)) onto the first \( r \) eigenvectors \( v_j \):
\begin{equation}
g \to P_S g = \sum_{j=1}^r \langle g, v_j \rangle v_j.
\end{equation}

The projection coefficients \( c_j \), which encode the contributions of each harmonic \( v_j \), are computed as:
\begin{equation}
c_j := \langle g, v_j \rangle = \langle \overline{W}, v_j \rangle = \begin{bmatrix}
 \langle \overline{W}_1, v_j \rangle \\
 \langle \overline{W}_2, v_j \rangle \\
 \vdots \\
 \langle \overline{W}_r, v_j \rangle
\end{bmatrix},
\end{equation}
\( \overline{W}_i \) is the \( i \)-th row of the matrix \( \overline{W} \). \( \langle \overline{W}_i, v_j \rangle = \sum_{k=1}^m (\overline{W}_i, \phi_k)v_j(\phi_k) \) represents projecting \( \overline{W}_i \) onto the \( j \)-th eigenvector \( v_j \). \( v_j(\phi_k) \) is the \( k \)-th component of the \( j \)-th eigenvector \( v_j \).

\textbf{Step 5: Latent harmonics.}
Calculating the latent harmonics at a new point \( \phi_{\text{new}} \):
\begin{equation}
V_j(\phi_{\text{new}}) = \sigma_j^{-1} \sum_{i=1}^m K_{\Phi}(\phi_{\text{new}}, \phi_i)v_j(\phi_i).
\end{equation}

\textbf{Step 6: Predicting in the Probabilistic Manifold Space}.
For a new point \( \phi_{\text{new}} \), the target function \( g(\phi_{\text{new}}) \) in the low-dimensional probabilistic manifold space is predicted as:
\begin{equation}
g(\phi_{\text{new}}) = \sum_{j=1}^r c_j V_j(\phi_{\text{new}}),
\end{equation}
where \( c_j \) are the projection coefficients computed in Step 4.
To predict the dynamics of the system, we set \( \phi_{\text{new}} = \Phi(x_{m+k-1}) \) to estimate \( \phi_{m+k} \), the latent representation of the data at the next time step:
\begin{equation}
\Phi(x_{m+k}) = g(\Phi(x_{m+k-1})).
\end{equation}

\textbf{Step 7: Reconstruction in the original space}.
Recovering the prediction in the original data space:
\begin{equation}
\mathbf{x}(m+k) = \psi(\Phi(x_{m+k})).
\end{equation}

\subsection{Nonlinear ROM via PMD} 

With the bidirectional mapping between the high-dimensional space and the low-dimensional probabilistic manifold, the reduced-order representation is given by:
\begin{equation} \label{eq:map}
u(t) = Q_r\hat{u}(t) + K\hat{X}(t) \in \mathbb{R}^n,
\end{equation}
where \( u(t) \) represents the high-dimensional velocity field governed by the Navier-Stokes equations, and \( Q_r \) and \( K \) are defined in equations \ref{eq:linear} and \ref{eq:kernal}. \( \hat{u}(t) \) denotes the low-dimensional representation of the linear features $\tilde{u}(t)$ (obtained via SVD-based projection), and \( \hat{X}(t) \) represents the low-dimensional representation of the nonlinear residual $X(t)$ defined in equation \ref{eq:error1}.

The ROM developed using the PMD method is highly flexible and applicable to various datasets without the need for extensive reconfiguration. During the reconstruction process, it can accurately capture the geometric and physical properties of the original data, ensuring reliable results.

In this way, the obtained low-dimensional data representation possesses both global stability (from the linear part) and precise retention of local geometric features (from the nonlinear part), making it better suited for modeling high-dimensional complex systems.

Applying equation \ref{eq:map} to the Navier-Stokes equations, we construct a ROM:
\begin{equation}
Q_r\nabla \cdot \hat{u}(t) + K\nabla \cdot \hat{X}(t) = 0,
\end{equation}
\begin{equation}
Q_r\overline{U}(t) + K\overline{X}(t) = -\nabla p + \nabla \cdot \tau,
\end{equation}
where the time-derivative and convective terms are given by:
\begin{equation}
\overline{U}(t) = \frac{\partial \hat{u}(t)}{\partial t} + \big( (Q_r\hat{u}(t) + K\hat{X}(t)) \cdot \nabla \big) \hat{u}(t) + \boldsymbol{f} \times \hat{u}(t),
\end{equation}
\begin{equation} 
\overline{X}(t) = \frac{\partial \hat{X}(t)}{\partial t} + \big( (Q_r\hat{u}(t) + K\hat{X}(t)) \cdot \nabla \big) \hat{X}(t) + \boldsymbol{f} \times \hat{X}(t),
\end{equation}
and the stress term is defined as:
\begin{equation} 
\nabla \cdot \tau = \nabla \cdot \mu \big( \nabla (Q_r\hat{u}(t) + K\hat{X}(t)) + (\nabla (Q_r\hat{u}(t) + K\hat{X}(t)))^T \big).
\end{equation}

\subsection{Dynamics prediction capability of the PMD}

PMD not only performs dimensionality reduction and reconstruction but also supports dynamic system's state predictions. It involves two steps: linear term prediction and then the nonlinear term.

To ensure the effectiveness of PMD in capturing the dominant dynamics, we analyze the energy captured by the basis functions obtained from both the linear and nonlinear feature extraction processes. As shown above, \( \sigma_i \) denote the singular values associated with the linear PMD basis (derived from linear features) and \( \lambda_i \) denote those corresponding to the nonlinear PMD basis (obtained from nonlinear feature analysis). For consistency, both the linear and nonlinear bases are truncated to a target dimension \( r \). The cumulative energy captured by the selected \( r \) modes is defined as:
\begin{equation}
E = \frac{\sum_{i=1}^{r} \sigma_i^2}{\sum_{i=1}^{m} \sigma_i^2} + \frac{\sum_{i=1}^{r} \lambda_i^2}{\sum_{i=1}^{m} \lambda_i^2}\cdot\frac{\sum_{i=r+1}^{m} \sigma_i^2}{\sum_{i=1}^{m} \sigma_i^2}.
\end{equation}
where \( m \) denotes the total number of available modes. The first term represents the fraction of the total linear energy captured by the top \( r \) modes, while the second term quantifies the contribution of the nonlinear modes relative to the residual linear energy. To maintain a balance between dimensionality reduction and information retention, we select the smallest \( r \) such that:
\begin{equation}
E \geq 0.95.
\end{equation}
This criterion ensures that at least 95\% of the system's variance is preserved in the reduced-order representation.

In comparison to the PMD method, the POD method typically employs the following energy capture formula for selecting the number of modes:

\begin{equation}
E_k = \frac{\sum_{i=1}^{k} \alpha_i^2}{\sum_{i=1}^{m} \alpha_i^2},
\end{equation}
where \( \alpha_i \) are the singular values associated with the POD basis and \( m \) is the total number of modes. The energy retained by the first \( k \) POD modes is computed by summing the squared singular values and normalizing by the total sum of squared singular values across all modes.

The key difference between the PMD and POD energy capture is that while POD solely relies on the singular values of the linear basis, PMD accounts for both the linear and nonlinear contributions. In PMD, the linear component captures the global dynamics, while the nonlinear component helps to model more complex, locally nonlinear behaviors. This leads to a more robust and accurate representation of the system, especially for highly nonlinear dynamics.

Thus, in practice, PMD provides a more comprehensive energy capture mechanism, enabling better handling of both linear and nonlinear dynamics, especially when the system exhibits significant nonlinearities that are not fully captured by POD. The improved accuracy and ability to maintain physical consistency in PMD are achieved by incorporating these nonlinear features in the modeling process.

\section{Convergence analysis of the PMD}
\label{sec:convergence}

A convergence analysis is crucial to validate the PMD method, ensuring that the sample matrix reliably approximates the underlying distribution matrix. This rigorous analysis enhances PMD's reliability and interpretability by offering a mathematically grounded framework for accuracy and computational efficiency.

\subsection{Convergence analysis of linearity in PMD}

This section examines the convergence of matrices capturing PMD's linear features, focusing on the decay of singular values and the convergence of sample covariance matrices to the true covariance matrix \( C \). PMD seeks an optimal low-dimensional representation by minimizing reconstruction error, demonstrating its efficiency and stability in handling high-dimensional data.

Given a data matrix \( U \in \mathbb{R}^{n \times m} \), where columns \( \mathbf{u}_j \) are independent and identically distributed (i.i.d) snapshots of the system state, the covariance matrix is:
\begin{equation}
C = \mathbb{E}[(\mathbf{U} - \mu)(\mathbf{U} - \mu)^T],
\end{equation}
where \( \mathbf{U} \) is a random vector, \( \mu \) is the true mean, and \( \mathbb{E}[\cdot] \) denotes the expectation operator.

To establish convergence, the sample covariance matrix is:
\begin{equation}
\hat{C}_m = \frac{1}{m} \sum_{i=1}^m (\mathbf{u}_i - \bar{\mathbf{u}})(\mathbf{u}_i - \bar{\mathbf{u}})^T,
\end{equation}
where \( \bar{\mathbf{u}} = \frac{1}{m} \sum_{i=1}^m \mathbf{u}_i \) is the sample mean. In the equation, \( m \to \infty \), \( \bar{\mathbf{u}} \to \mu \), and \( \hat{C}_m \to C \). Furthermore, the eigenvalues and eigenvectors of \( \hat{C}_m \) converge to those of \( C \), ensuring the precision of the linear features of the PMD.

\subsubsection{Convergence of the sample covariance matrix}
We analyze the convergence rate of the sample covariance matrix \( \hat{C}_m \) to the true covariance matrix \( C \), showing that \( \mathbb{E}[\| \hat{C}_m - C \|] = \mathcal{O}(\frac{1}{\sqrt{m}}) \), where \( \| \cdot \| \) is the spectral norm. The convergence analysis uses the law of large numbers (LLN)\cite{hsu1947complete}, central limit theorem (CLT)\cite{rosenblatt1956central}, and matrix perturbation theory\cite{baumgartel1984analytic}. The sample covariance matrix is defined as:
\begin{equation}
\hat{C}_m = \frac{1}{m} \sum_{i=1}^m (\mathbf{u}_i - \mu)(\mathbf{u}_i - \mu)^T,
\end{equation}
where \( \mathbf{u}_i \) is i.i.d samples with true covariance \( C \). By the LLN:
\begin{equation}
\hat{C}_m \stackrel{a.s.}{\to} C \quad \text{as } m \to \infty.
\end{equation}
To determine the convergence rate, we define the error matrix:
\begin{equation}\label{eq:error2}
E_m = \hat{C}_m - C.
\end{equation}
Each element \( (E_m)_{ab} \) is a sample mean of i.i.d random variables, and by the CLT:
\begin{equation}
\operatorname{Var}\left( (E_m)_{ab} \right) = \frac{1}{m} \operatorname{Var}\left( (\mathbf{u}_i - \mu)_a (\mathbf{u}_i - \mu)_b \right).
\end{equation}
Thus, the standard deviation of each entry scales as \( \mathcal{O}(\frac{1}{\sqrt{m}}) \).

Using matrix concentration inequalities, the spectral norm \( \| E_m \| \), which measures the largest singular value of \( E_m \), also scales as \( \mathcal{O}(\frac{1}{\sqrt{m}}) \). Therefore:
\begin{equation}\label{eq:order}
\mathbb{E}[\| \hat{C}_m - C \|] = \mathcal{O}(\frac{1}{\sqrt{m}}).
\end{equation}
This establishes that the sample covariance matrix converges to the true covariance matrix with a rate of \( \mathcal{O}(\frac{1}{\sqrt{m}}) \) in the spectral norm.

\subsubsection{Convergence of eigenvalues and eigenvectors}

Given that \(\hat{C}_m\) converges to \(C\) at a rate of \(\mathcal{O}(\frac{1}{\sqrt{m}})\), we analyze the convergence of its eigenvalues and eigenvectors to those of \(C\). Let \(\lambda_1 \geq \lambda_2 \geq \dots \geq \lambda_n\) and \(v_1, v_2, \dots, v_n\) denote the eigenvalues and eigenvectors of \(C\), and \(\hat{\lambda}_1 \geq \hat{\lambda}_2 \geq \dots \geq \hat{\lambda}_n\) and \(\hat{v}_1, \hat{v}_2, \dots, \hat{v}_n\) those of \(\hat{C}_m\). The linear bases in equation \ref{eq:linear} correspond to the eigenvectors of \(\hat{C}_m\).

The convergence of eigenvectors depends on both sample size \(m\) and eigenvalue gaps \(\lambda_i - \lambda_{i+1}\). Using matrix perturbation theory \cite{stewart1998perturbation} and Davis-Kahan theorem \cite{yu2015useful}, we estimate:
\begin{equation}\label{eq:vectors}
\left\|\hat{v}_i - v_i\right\| \sim \mathcal{O}\left(\frac{1}{\sqrt{m}} \cdot \frac{1}{\lambda_i - \lambda_{i+1}}\right),
\end{equation}
indicating slower convergence when eigenvalues are close. Similarly, eigenvalue deviations satisfy:
\begin{equation}
|\hat{\lambda}_i - \lambda_i| \leq \|E\|_2,
\end{equation}
where \(E = \hat{C}_m - C\) and \(\|E\|_2 \sim \mathcal{O}(\frac{1}{\sqrt{m}})\), as shown in Weyl's inequality \cite{fan1949theorem}.
Key convergence results include the following:
\begin{itemize}
    \item \textbf{Matrix convergence}: \(\mathbb{E}[|\hat{C}_m - C|] = \mathcal{O}(\frac{1}{\sqrt{m}})\).
    \item \textbf{Eigenvalue convergence}: \(\mathbb{E}[|\hat{\lambda}_i - \lambda_i|] = \mathcal{O}(\frac{1}{\sqrt{m}})\).
    \item \textbf{Eigenvector convergence}: If \(\delta = \lambda_i - \lambda_{i+1} > 0\), then \(\mathbb{E}[\|\hat{v}_i - v_i\|] = \mathcal{O}\left(\frac{1}{\delta\sqrt{m}}\right)\).
\end{itemize}
As \(m \to \infty\), \(\hat{\lambda}_i \to \lambda_i\) and \(\hat{v}_i \to v_i\), they ensure that the PMD method accurately captures the linear and nonlinear dynamics.

\subsection{Convergence analysis of the non-linearity in PMD}
The PMD method, rooted in spectral graph theory and kernel methods, addresses the convergence of both linear and nonlinear features. This analysis focuses on the convergence of kernel density estimation, spectral decomposition, and the effect of manifold structure on convergence behavior.

Consider the dataset of equation \ref{eq:error1} with \( m \) data points \( \{ \mathbf{x}_1, \mathbf{x}_2, \dots, \mathbf{x}_m \} \subset \mathbb{R}^n \), sampled from a distribution \( p(x) \) on a low-dimensional manifold \( \mathcal{M} \). The goal of PMD is to approximate the manifold structure and underlying probability distribution while maintaining computational efficiency and accuracy.

\subsubsection{Convergence of the probability distribution}

PMD constructs a sample similarity matrix $\hat{K}$ and compares it with the true similarity matrix $K$, derived from the probability distribution $p(x)$. As the sample size $m$ increases, $\hat{K}$ converges to $K$, revealing the structure of the manifold. This convergence relies on kernel density estimation (KDE). Given data points$\{x_i\}_{i=1}^m$, and a kernel function \( k_{\epsilon}(x, y) \). Assume we use a standard Gaussian kernel function, $x,y$ are arbitrary points on the manifold:
\begin{equation}
k_{\epsilon}(x,y) = \frac{1}{(2\pi\epsilon^2)^{d/2}} \exp\left(-\frac{\|x-y\|^2}{2\epsilon^2}\right),
\end{equation}
and the KDE is:
\begin{equation}
\hat{p}(x) = \frac{1}{m\epsilon^d} \sum_{i=1}^{m} k_{\epsilon}(x, x_i),
\end{equation}
where $\epsilon$ is the bandwidth and $d$ is the intrinsic dimension, $m \rightarrow \infty$, $\hat{p}(x) \rightarrow p(x)$ under the conditions $\epsilon \rightarrow 0$ and $m\epsilon^d \rightarrow \infty$.

To analyze the convergence, the expectation of $\hat{p}(x)$ is:
\begin{equation}
E[\hat{p}(x)] = \int p(y) k_{\epsilon}(x, y) \, dy.
\end{equation}
We can expand \( p(y) \) in a Taylor series around \( x \):
\begin{equation}
p(y) \approx p(x) + (y - x)^T \nabla p(x) + \frac{1}{2}(y - x)^T H(p(x))(y - x),
\end{equation}
where $\nabla p(x)$ is the gradient, and $H(p(x))$ is the Hessian matrix.
\\
\textbf{Key Terms in the Expansion:}
\begin{enumerate}
    \item \textbf{Constant Term:} $\int p(x) k_{\epsilon}(x, y) \, dy = p(x).$
    \item \textbf{Gradient Term:} $\int (y - x)^T \nabla p(x) k_{\epsilon}(x, y) \, dy = 0$ (symmetry of $k_{\epsilon}$).
    \item \textbf{Hessian Term:} Using the kernel's second moment property:
    \begin{equation}
    \int (y_i - x_i)(y_j - x_j) k_{\epsilon}(x, y) \, dy = \epsilon^2 \delta_{ij},
    \end{equation}
    the Hessian term becomes:
    \begin{equation}
    \frac{1}{2} \epsilon^2 \text{tr}(H(p(x))),
    \end{equation}
    where $\text{tr}(H(p(x)))$ is the trace of the Hessian.
\end{enumerate}

Thus, the bias is:
\begin{equation}
\text{Bias}[\hat{p}(x)] = E[\hat{p}_{\epsilon}(x)] - p(x) \sim \mathcal{O}(\epsilon^2).
\end{equation}
The variance of kernel density estimation is caused by the deviation of each independent sample from the estimated value. Its variance is:
\begin{equation}
\text{Var}[\hat{p}(x)] \sim \frac{1}{m\epsilon^d},
\end{equation}
Since the variance of each kernel function is \( \mathcal{O}(1) \) and there are \( m \) samples, the variance decreases at a rate of \( \mathcal{O}(\frac{1}{m\epsilon^d}) \) as \( m \rightarrow \infty \).
\\
\textbf{Mean Squared Error (MSE):}
Combining bias and variance:
\begin{equation}
\text{MSE} = \text{Bias}^2 + \text{Var} \sim \mathcal{O}(\epsilon^4) + \mathcal{O}\left(\frac{1}{m\epsilon^d}\right).
\end{equation}
To minimize the MSE, balance the terms:
\begin{equation}
\epsilon^4 \sim \frac{1}{m\epsilon^d}.
\end{equation}
Solving for $\epsilon$:
\begin{equation}
\epsilon \sim m^{-\frac{1}{d+4}}.
\end{equation}
Substituting this $\epsilon$ into the MSE gives the convergence rate:
\begin{equation}
\text{MSE} \sim \mathcal{O}\left(m^{-\frac{4}{d+4}}\right).
\end{equation}
This rate demonstrates that KDE accurately approximates the true density $p(x)$ as $m$ grows.

\subsubsection{Convergence of the transition matrix}
The PMD method leverages spectral decomposition of a normalized similarity matrix \( P_{\epsilon} \), defined as:
\begin{equation}  
P_{\epsilon} = D^{-1}\hat{\mathbf{K}}, \quad D_{ii} = \sum_j \hat{\mathbf{K}}_{ij}.
\end{equation}

As \( m \to \infty \), \( P_\epsilon \) converges to a continuous operator \( P_\infty \):
\begin{equation} 
P_\infty f(x) = \frac{\int_\mathcal{M} k_\epsilon(x, y) f(y) p(y) \, d\mu(y)}{\int_\mathcal{M} k_\epsilon(x, y) p(y) \, d\mu(y)},
\end{equation}
where \( \mathcal{M} \) is the manifold, \( p(y) \) is the sample density, $f(y)$ is the test function, and \( d\mu(y) \) the measure. For finite samples, integrals are approximated as:
\begin{equation} 
\int_\mathcal{M} f(y) p(y) \, d\mu(y) \approx \frac{1}{m} \sum_{j=1}^m f(x_j).
\end{equation} 
\\
\textbf{Law of Large Numbers.}  
For i.i.d samples \( x_1, \dots, x_m \), the sample mean:
\begin{equation} 
\hat{I} = \frac{1}{m} \sum_{i=1}^m f(x_i),
\end{equation} 
converges to the true integral \( I = \int_\mathcal{M} f(y) p(y) \, d\mu(y) \) almost surely as \( m \to \infty \):
\begin{equation} 
\hat{I} \xrightarrow{\text{a.s.}} I.
\end{equation} 
\\
\textbf{Variance and Error Bound.}  
For bounded \( f(y) \), the variance of the sample mean is:
\begin{equation} 
\text{Var}(\hat{I}) = \frac{\sigma^2}{m}, \quad \sigma^2 = \int_\mathcal{M} \left(f(y) - I\right)^2 p(y) \, d\mu(y).
\end{equation} 
Using Hoeffding's inequality, the error probability satisfies:
\begin{equation} 
\mathbb{P}\left(\left|\hat{I} - I\right| \geq \delta\right) \leq 2 \exp\left(-\frac{2m\delta^2}{M^2}\right),
\end{equation} 
where \( \|f(y)\|_\infty \leq M \). Thus, the error decreases as:
\begin{equation} 
\left|\hat{I} - I\right| \leq \frac{C}{\sqrt{m}},
\end{equation} 
where \( C \) depends on \( \sigma^2 \) and \( M \).

For PMD, kernel integrals are approximated as:
\begin{equation} 
\int_\mathcal{M} k_\epsilon(x, y) p(y) \, d\mu(y) \approx \frac{1}{m} \sum_{j=1}^m k_\epsilon(x, x_j).
\end{equation} 
The error in \( P_\epsilon \) and \( P_\infty \) satisfies:
\begin{equation}
\|P_\epsilon - P_\infty\| \leq \frac{C'}{\sqrt{m}},
\end{equation}
where \( \|\cdot\| \) is the operator norm.

\subsubsection{Convergence of eigenvalues and eigenvectors}

Suppose $\lambda_i^{(m)}$ is the eigenvalue of $P_\epsilon$, $\psi_i^{(m)}$ is the eigenvector corresponding to it, and $\lambda_{i}^{(\infty)}$ is the corresponding eigenvalue of $P_\infty$, $\psi_i^{\infty}$ is the eigenvector corresponding to it. By matrix perturbation theory and the analysis shown before, we have:
\begin{equation} \label{eq:verror}
|\lambda_i^{(m)} - \lambda_{i}^{(\infty)}| \leq \|P_\epsilon - P_\infty\| \leq \frac{C_1}{\sqrt{m}}.
\end{equation}

Similarly, for eigenvectors $\psi(x)$, the error between the finite-sample eigenvector and the theoretical eigenfunction satisfies:
\begin{equation}
\|\psi_i^{(m)}(x_j) - \psi_i^{\infty}(x_j)\| \leq \frac{C_2}{\sqrt{m}},
\end{equation}
where $C_2$ depends on the smoothness of the kernel function, the regularity of eigenvectors, the gap between the eigenvalues, and normalization conditions.

\subsubsection{Convergence of the probabilistic manifold distance}

The probabilistic manifold distance \( D_t(x, y) \) quantifies the similarity between two points \( x \) and \( y \) on the probabilistic manifold after \( t \) steps:
\begin{equation}
D_t^{(\infty)}(x, y) = \sqrt{\sum_{i=1}^\infty {\lambda_i^{(\infty)}}^{2t} \left(\psi_i^{(\infty)}(x) - \psi_i^{(\infty)}(y)\right)^2},
\end{equation}
where \( \lambda_i \) are the eigenvalues of the PMD operator, \( \psi_i(x) \) are the eigenvectors, and \( t \) is the time step.
For finite samples, it is approximated as:
\begin{equation}
D_t^{(m)}(x, y) = \sqrt{\sum_{i=1}^r \left(\lambda_i^{(m)}\right)^{2t} \left(\psi_i^{(m)}(x) - \psi_i^{(m)}(y)\right)^2}.
\end{equation}

The error between \( D_t^{(m)}(x, y) \) and \( D_t^{(\infty)}(x, y) \) is:
\begin{equation}
\Delta D_t = \left| (D_t^{(m)})^2(x, y) - (D_t^{(\infty)})^2(x, y) \right| = O\left(\frac{1}{\sqrt{m}}\right),
\end{equation}
which can be expressed as:
\begin{equation}
\Delta D_t \leq D_1 + D_2 + D_3.
\end{equation}

The components of \( \Delta D_t \) are as follows:
\\
\textbf{1. Error from eigenvalues}
\begin{equation}
D_1 = \sum_{i=1}^r \left| \left(\lambda_i^{(m)}\right)^{2t} - \left(\lambda_i^{(\infty)}\right)^{2t} \right| \cdot \left(\psi_i^{(\infty)}(x) - \psi_i^{(\infty)}(y)\right)^2.
\end{equation}
Using eigenvalue convergence:
\begin{equation}
|\lambda_i^{(m)} - \lambda_i^{(\infty)}| \leq \frac{C_1}{\sqrt{m}},
\end{equation}
the propagated error is:
\begin{equation}
\left| \left(\lambda_i^{(m)}\right)^{2t} - \left(\lambda_i^{(\infty)}\right)^{2t} \right| \leq 2t \cdot \left(\lambda_i^{(\infty)}\right)^{2t-1} \cdot \frac{C_1}{\sqrt{m}}.
\end{equation}
Thus:
\begin{equation}
D_1 \leq \frac{C_1^{'}}{\sqrt{m}}.
\end{equation}
\\
\textbf{2. Error from eigenvectors}
\begin{equation}
D_2 = \sum_{i=1}^r \left(\lambda_i^{(m)}\right)^{2t} \cdot \left| (\psi_i^{(m)}(x) - \psi_i^{(m)}(y))^2 - (\psi_i^{(\infty)}(x) - \psi_i^{(\infty)}(y))^2 \right|.
\end{equation}
The eigenvector error satisfies:
\begin{equation}
\|\psi_i^{(m)} - \psi_i^{(\infty)}\| \leq \frac{C_2}{\sqrt{m}}.
\end{equation}
Thus:
\begin{equation}
D_2 \leq \frac{C_2^{'}}{\sqrt{m}}.
\end{equation}
\\
\textbf{3. Truncation error}.
In practice, only the first \( r \) terms of the series are computed, neglecting higher-order terms. The truncation error is:
\begin{equation}
D_3 = \sum_{i=r+1}^\infty \left(\lambda_i^{(\infty)}\right)^{2t} \left(\psi_i^{(\infty)}(x) - \psi_i^{(\infty)}(y)\right)^2,
\end{equation}
which can be made arbitrarily small if the eigenvalues decay rapidly (\( \lambda_i^{(\infty)} \to 0 \) as \( i \to \infty \)) by choosing a sufficiently large \( r \).

To confirm the exponential decay of eigenvalues, we connect PMD with the geometric structure of the manifold. The heat equation describes the process:
\begin{equation}
\frac{\partial u(x, t)}{\partial t} = \Delta u(x, t),
\end{equation}
where \( u(x, t) \) represents the diffusion state at position \( x \) and time \( t \), and \( \Delta \) is the Laplace-Beltrami operator. This equation highlights how the manifold's geometric properties, such as curvature, influence the diffusion process.\\
\\
\textbf{Spectral Decomposition of the Laplace-Beltrami Operator}

The spectral theory of the Laplace-Beltrami operator resolves any diffusion process on a manifold via eigenfunctions and eigenvalues:
\begin{equation}
\Delta \phi_i(x) = -\mu_i \phi_i(x),
\end{equation}
where \( \phi_i(x) \) is the \( i \)-th eigenfunction, and \( \mu_i \) is its eigenvalue.

The heat equation solution is:
\begin{equation}
u(x, t) = e^{-t\Delta} u(x, 0),
\end{equation}
which can be expanded spectrally:
\begin{equation}
u(x, t) = \sum_{i=1}^\infty e^{-\mu_i t} \phi_i(x) \langle u(x, 0), \phi_i(x) \rangle.
\end{equation}
As \( t \) increases, eigenmodes decay at rates determined by \( e^{-\mu_i t} \), with larger \( \mu_i \) leading to faster decay.  

The eigenvalues \( \lambda^{(\infty)}_i \) of the transition matrix \( P_\infty \) exhibit similar decay behavior. The temporal evolution of the PMD process can be written as:
\begin{equation}
p(t) = P_\infty^t p(0).
\end{equation}
With eigenvalue decomposition \( P_\infty = \sum_{i=1}^m \lambda^{(\infty)}_i \psi^{(\infty)}_i {\psi^{(\infty)}_i}^T \), the evolution becomes:
\begin{equation}
p(t) = \sum_{i=1}^m {\lambda^{(\infty)}_i}^t \psi^{(\infty)}_i {\psi_i^{(\infty)}}^T p(0).
\end{equation}

For eigenvalues \( \lambda^{(\infty)}_i \) close to 1:
\begin{equation}
{\lambda_i^{(\infty)}}^t \approx \exp(-\mu_i t),
\end{equation}
where \( \mu_i \) relates to the manifold's diffusion rate. Thus, eigenvalues decay exponentially:
\begin{equation}
\lambda^{(\infty)}_i(t) \sim \exp(-\mu_i t).
\end{equation}

This result reflects the local geometric structure of the manifold and confirms that the truncation error \( D_3 \) approaches zero as the number of retained eigenmodes increases:
\begin{equation}
D_3 \to 0.
\end{equation}
\\
\textbf{4. Total Error}. Combining all three sources of error, the total error in diffusion distance is:
\begin{equation}
\Delta D_t \leq D_1 + D_2 + D_3.
\end{equation}

Substituting the bounds:
\begin{equation}
\left| D_t^{(n)}(x, y) - D_t^{(\infty)}(x, y) \right| \leq \frac{C_0}{\sqrt{m}},
\end{equation}
where \( C_0 \) is a constant that depends on the eigenvalue decay, eigenfunction regularity, and manifold properties.

\begin{itemize}
    \item \textbf{Kernel Density Estimation}: The KDE converges to the true density at a rate of \( \mathcal{O}\left(m^{-\frac{4}{d+4}}\right) \).
    \item \textbf{Transition Matrix}: \( P_\epsilon \) converge to the continuous operator at a rate of \( \mathcal{O}(\frac{1}{\sqrt{m}}) \).
    \item \textbf{Eigenvalues and Eigenvectors}: Eigenvalues and eigenvectors of \( P_\epsilon \) converge to those of the continuous operator at a rate of \( \mathcal{O}(\frac{1}{\sqrt{m}}) \).
    \item \textbf{Probabilistic Manifold Distance}: The probabilistic manifold distance converges to the geodesic distance on \( \mathcal{M} \) as \( m \to \infty \) ar a rate of  \( \mathcal{O}(\frac{1}{\sqrt{m}}) \).
\end{itemize}

These convergence results highlight how PMD effectively captures the nonlinear probabilistic and geometric structure of the data. The convergence rates are shaped by factors such as sample size, manifold dimensionality, and the properties of the kernel bandwidth, enabling PMD to model complex data relationships with high fidelity.

\section{Numerical experiments} \label{sec:experiments}
A demonstration of the use of the PMD based ROM is presented in this section. This is based on solving the Navier–Stokes equations for lock exchange and flow past a cylinder problems. The PMD is implemented under the framework of Fluidity and compared against ROM based on POD+LSTM and the full high-fidelity model. 
\subsection{Lock exchange problem}
The lock exchange problem, or gravity current problem, is a benchmark in fluid dynamics that models the interaction of two fluids with different densities. This scenario typically involves a rectangular tank divided by a barrier, with denser fluid (e.g., saltwater) on one side and lighter fluid (e.g., freshwater) on the other. When the barrier is removed, gravity-driven mixing occurs, producing complex flow patterns such as vortices and waves.

The governing Navier-Stokes equations make this problem applicable to studies of oceanic mixing, atmospheric flows, and industrial processes. In the numerical example, a two-dimensional setup follows \cite{shin2004gravity}, with the denser fluid initially on the left. The dimensionless domain (\(0.8 \times 0.1\)) contains 1491 grid nodes.

The simulation has 40 seconds and 160 snapshots at uniform time intervals \(\Delta t = 0.25\). These snapshots include velocity components and pressure. 120 snapshots(30s) were used to construct the ROM, and the remaining 40 snapshots are used to test the predictive capabilities of the ROMs. 

\begin{figure}[tbhp] 
\centering
\begin{tabular}{cc}

\begin{minipage}{0.45\linewidth} 
\includegraphics[width = \linewidth,angle=0,clip=true]{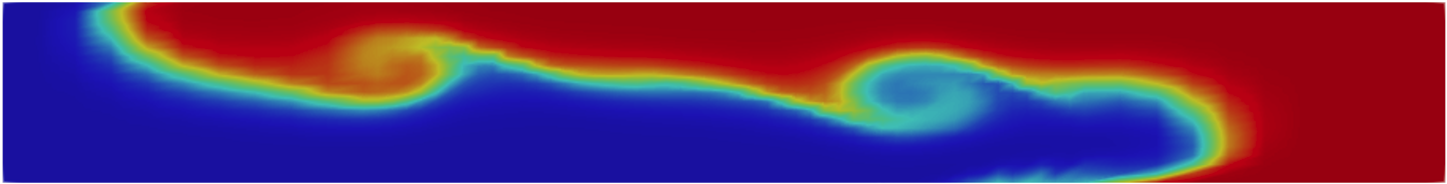}
\end{minipage}
&
\begin{minipage}{0.45\linewidth} 
\includegraphics[width = \linewidth,angle=0,clip=true]{lockexchangefullmodelt=1.png}
\end{minipage}\\
(a) {\small Full model, $t = 25$s}&
(b) {\small Full model, $t = 25$s}\\

\begin{minipage}{0.45\linewidth} 
\includegraphics[width = \linewidth,angle=0,clip=true]{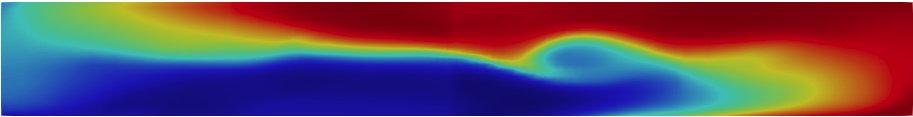}
\end{minipage}
&
\begin{minipage}{0.45\linewidth} 
\includegraphics[width = \linewidth,angle=0,clip=true]{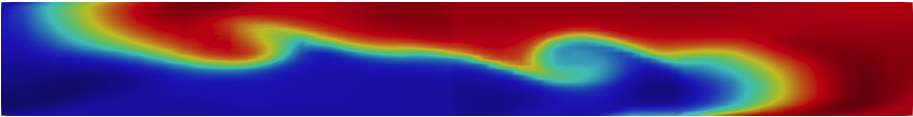}
\end{minipage}\\
(c) {\small POD+LSTM, $r = 2, t = 25$s}&
(d) {\small POD+LSTM, $r = 6, t = 25$s}\\

\begin{minipage}{0.45\linewidth} 
\includegraphics[width = \linewidth,angle=0,clip=true]{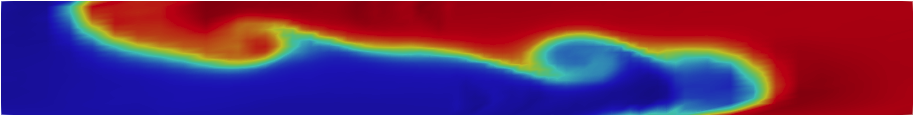}
\end{minipage}
&
\begin{minipage}{0.45\linewidth} 
\includegraphics[width = \linewidth,angle=0,clip=true]{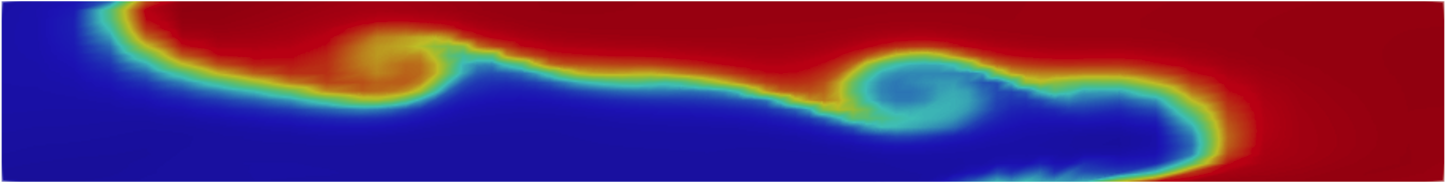}
\end{minipage}\\
(e) {\small PMD, $r = 2, t = 25$s}&
(f) {\small PMD, $r = 6, t = 25$s}\\ 

\begin{minipage}{0.45\linewidth} 
\includegraphics[width = \linewidth,angle=0,clip=true]{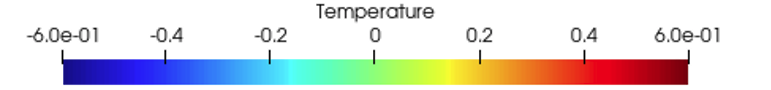}
\end{minipage}
&
\begin{minipage}{0.45\linewidth} 
\includegraphics[width = \linewidth,angle=0,clip=true]{lockscale1.png}
\end{minipage}\\

\end{tabular}
\caption{Lock exchange test case: temperature values obtained from the high fidelity full model, POD+LSTM and PMD at $t = 25$ using 2 ($r=2$) and 6 basis functions ($r=6$) respectively.}\label{fg:lockconstruct}
\end{figure}


\begin{figure}[tbhp]
\centering
\begin{tabular}{cc}
\begin{minipage}{0.45\linewidth}
\includegraphics[width=\linewidth]{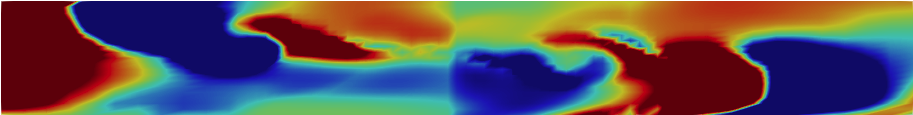}
\end{minipage} &
\begin{minipage}{0.45\linewidth}
\includegraphics[width=\linewidth]{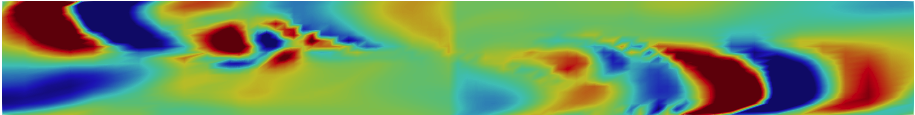}
\end{minipage} \\
(a) {\small POD, $r = 2, t = 25$s} & (b) {\small POD, $r = 6, t = 25$s} \\

\begin{minipage}{0.45\linewidth}
\includegraphics[width=\linewidth]{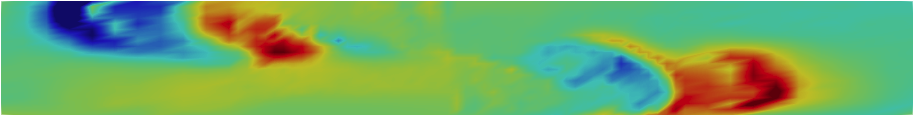}
\end{minipage} &
\begin{minipage}{0.45\linewidth}
\includegraphics[width=\linewidth]{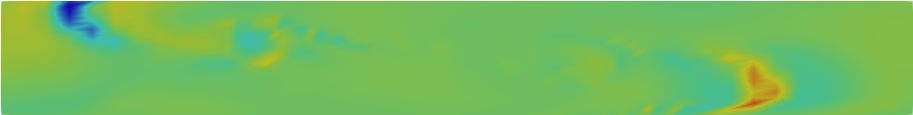}
\end{minipage} \\
(c) {\small PMD, $r = 2, t = 25$s} & (d) {\small PMD, $r = 6, t = 25$s} \\

\begin{minipage}{0.45\linewidth} 
\includegraphics[width = \linewidth,angle=0,clip=true]{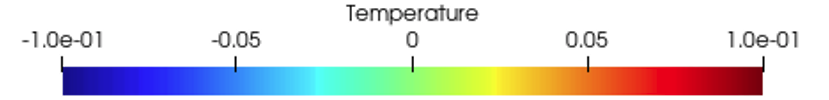}
\end{minipage}
&
\begin{minipage}{0.45\linewidth} 
\includegraphics[width = \linewidth,angle=0,clip=true]{lockscale2.png}
\end{minipage}\\
\end{tabular}
\caption{Lock exchange test case: errors of PMD+LSTM and PMD using 2 and 6 basis functions.}\label{fg:lockerror1}
\end{figure}

\ref{fg:lockconstruct} and \ref{fg:lockerror1} compare the temperature results of the lock exchange test case obtained from the high fidelity full model, PMD and ROM based on POD and LSTM using 4 and 8 basis functions. The error analysis is also conducted. As shown in the figures, POD method struggles to capture nonlinear dynamics using 4 basis functions while PMD performs better and captures most of the details. The accuracy of POD can be increased by increasing the number of basis functions, which is shown in \ref{fg:lockconstruct} (d). This is similar to PMD method. As shown in  \ref{fg:lockerror1}, the errors of PMD using 8 basis functions are smaller than those of PMD using 4 basis functions and POD using 8 basis functions. 
These results highlight PMD performs better than POD based ROM using the same number of basis functions.   



\begin{figure}[tbhp]
\centering
\begin{tabular}{cc}
\begin{minipage}{0.45\linewidth}
\includegraphics[width=\linewidth]{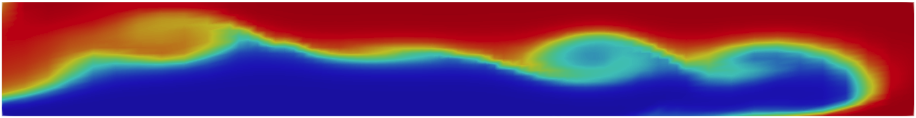}
\end{minipage} &
\begin{minipage}{0.45\linewidth}
\includegraphics[width=\linewidth]{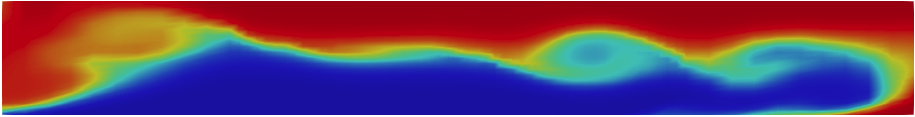}
\end{minipage} \\
(a) {\small Full model, $t = 33$s} & (b) {\small Full model, $t = 36$s} \\

\begin{minipage}{0.45\linewidth}
\includegraphics[width=\linewidth]{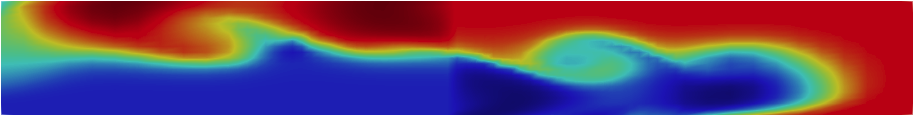}
\end{minipage} &
\begin{minipage}{0.45\linewidth}
\includegraphics[width=\linewidth]{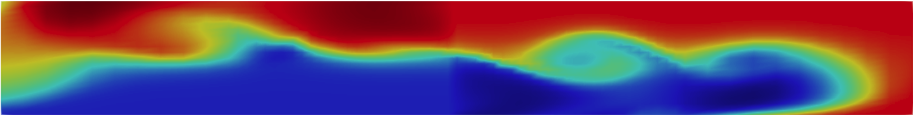}
\end{minipage} \\
(c) {\small POD+LSTM, $r = 4, t = 33$s} & (d) {\small POD+LSTM, $r = 4, t = 36$s} \\

\begin{minipage}{0.45\linewidth}
\includegraphics[width=\linewidth]{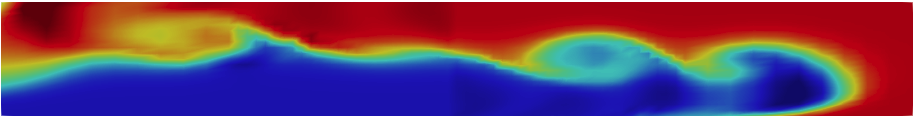}
\end{minipage} &
\begin{minipage}{0.45\linewidth}
\includegraphics[width=\linewidth]{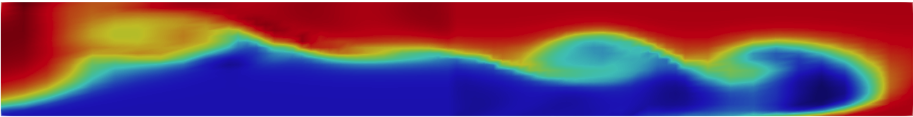}
\end{minipage} \\
(e) {\small POD+LSTM, $r = 8, t = 33$s} & (f) {\small POD+LSTM, $r = 8, t = 36$s} \\

\begin{minipage}{0.45\linewidth}
\includegraphics[width=\linewidth]{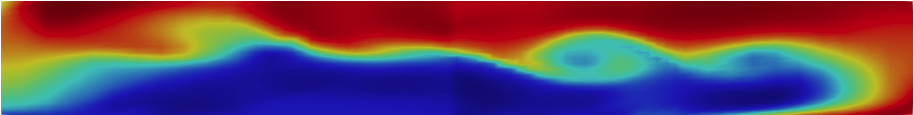}
\end{minipage} &
\begin{minipage}{0.45\linewidth}
\includegraphics[width=\linewidth]{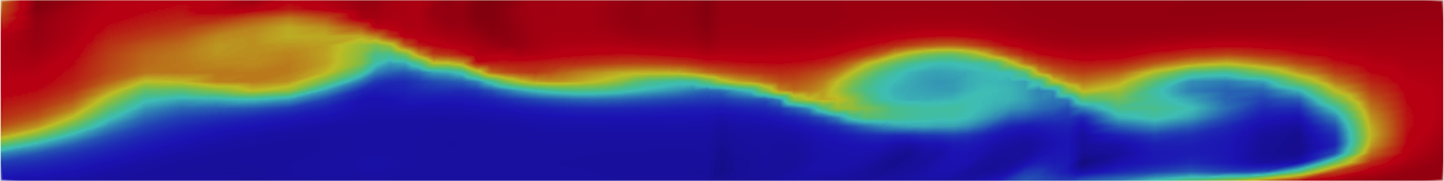}
\end{minipage} \\
(g) {\small PMD, $r = 4, t = 33$s} & (h) {\small PMD, $r = 4, t = 36$s} \\

\begin{minipage}{0.45\linewidth}
\includegraphics[width=\linewidth]{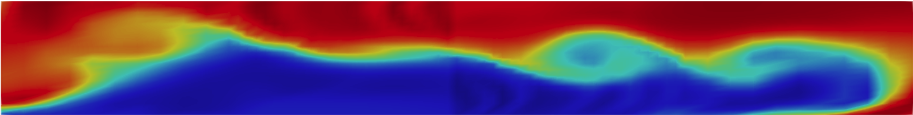}
\end{minipage} &
\begin{minipage}{0.45\linewidth}
\includegraphics[width=\linewidth]{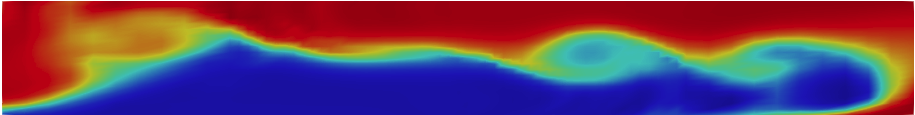}
\end{minipage} \\
(i) {\small PMD, $r = 8, t = 33$s} & (j) {\small PMD, $r = 8, t = 36$s} \\

\begin{minipage}{0.45\linewidth} 
\includegraphics[width = \linewidth,angle=0,clip=true]{lockscale1.png}
\end{minipage}
&
\begin{minipage}{0.45\linewidth} 
\includegraphics[width = \linewidth,angle=0,clip=true]{lockscale1.png}
\end{minipage}\\
 
\end{tabular}
\caption{Lock exchange test case: temperature solutions at $t = 33s$ and $36s$ using POD+LSTM and PMD using 4 and 8 basis functions respectively.}
\label{fg:lockpredict}
\end{figure}

\begin{figure}[tbhp]
\centering
\begin{tabular}{cc}
\begin{minipage}{0.45\linewidth}
\includegraphics[width=\linewidth]{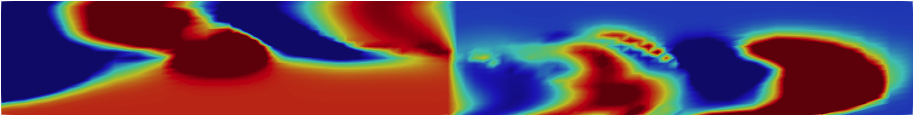}
\end{minipage} &
\begin{minipage}{0.45\linewidth}
\includegraphics[width=\linewidth]{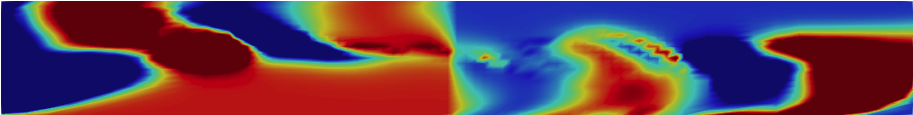}
\end{minipage} \\
(a) {\small POD+LSTM, $r = 4, t = 33s$} & (b) {\small POD+LSTM, $r = 4, t = 36s$} \\

\begin{minipage}{0.45\linewidth}
\includegraphics[width=\linewidth]{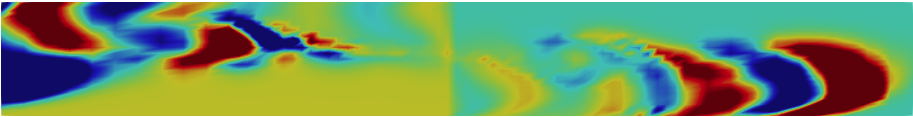}
\end{minipage} &
\begin{minipage}{0.45\linewidth}
\includegraphics[width=\linewidth]{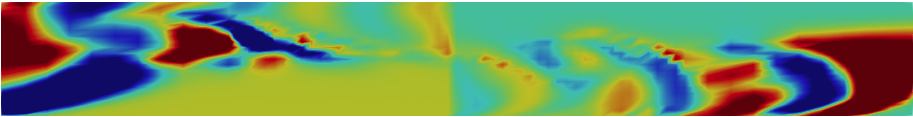}
\end{minipage} \\
(c) {\small POD+LSTM, $r = 8, t = 33s$} & (d) {\small POD+LSTM, $r = 8, t = 36s$} \\

\begin{minipage}{0.45\linewidth}
\includegraphics[width=\linewidth]{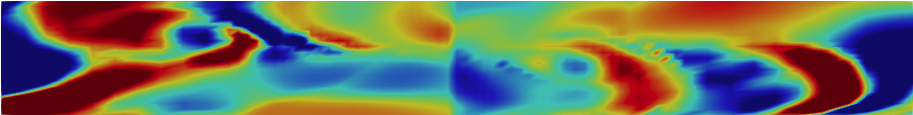}
\end{minipage} &
\begin{minipage}{0.45\linewidth}
\includegraphics[width=\linewidth]{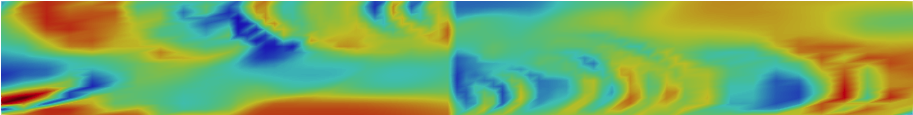}
\end{minipage} \\
(e) {\small PMD, $r = 4, t = 33s$} & (f) {\small PMD, $r = 4, t = 36s$} \\

\begin{minipage}{0.45\linewidth}
\includegraphics[width=\linewidth]{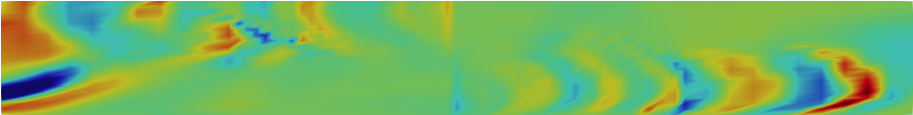}
\end{minipage} &
\begin{minipage}{0.45\linewidth}
\includegraphics[width=\linewidth]{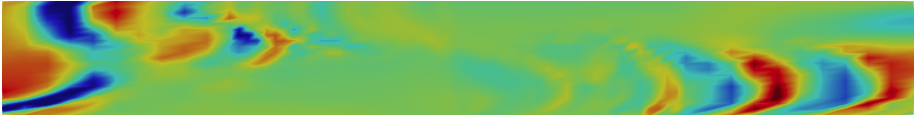}
\end{minipage} \\
(g) {\small PMD, $r = 8, t = 33s$} & (h) {\small PMD, $r = 8, t = 36s$} \\

\begin{minipage}{0.45\linewidth} 
\includegraphics[width = \linewidth,angle=0,clip=true]{lockscale2.png}
\end{minipage}
&
\begin{minipage}{0.45\linewidth} 
\includegraphics[width = \linewidth,angle=0,clip=true]{lockscale2.png}
\end{minipage}\\

\end{tabular}
\caption{Lock exchange test case: errors of POD+LSTM and PMD at $t = 33s$ and $36s$ using 4 and 8 basis functions respectively.}
\label{fg:lockerror2}
\end{figure}

In order to demonstrate the predictive capabilities of the PMD,  two unseen time levels' results ($t=33$s and $t=36$s) outside the range of time steps that were used to construct the ROM are given in \ref{fg:lockpredict}. The figure shows that the PMD predicts better than POD based ROM using the same number of basis functions. In order to see the accuracy differences between the POD based ROM and PMD, the errors of POD+LSTM and PMD are given in \ref{fg:lockerror2}. As shown in this figure, both POD+LSTM and PMD can increase the accuracy via using a larger number of basis functions. In addition, PMD has less errors than those of POD+LSTM using the same number of basis functions or modes. 


\subsection{Flow past a cylinder}

The \textit{flow around a cylinder} is a classic CFD benchmark problem for studying flow separation, vortex shedding, and wake dynamics. Fluid behavior is examined as it flows past a stationary cylinder at varying Reynolds numbers (\(Re\)), defined as:  
$Re = \frac{\rho UD}{\mu}$, where \(U\) is the inflow velocity, \(\rho\) the fluid density, \(D\) the cylinder diameter, and \(\mu\) the dynamic viscosity.
The PMD method was applied to simulate and analyze this problem, effectively capturing non-linearity.  

\begin{figure}[tbhp]
\centering
\begin{tabular}{cc}
\begin{minipage}{0.45\linewidth}
\includegraphics[width=\linewidth]{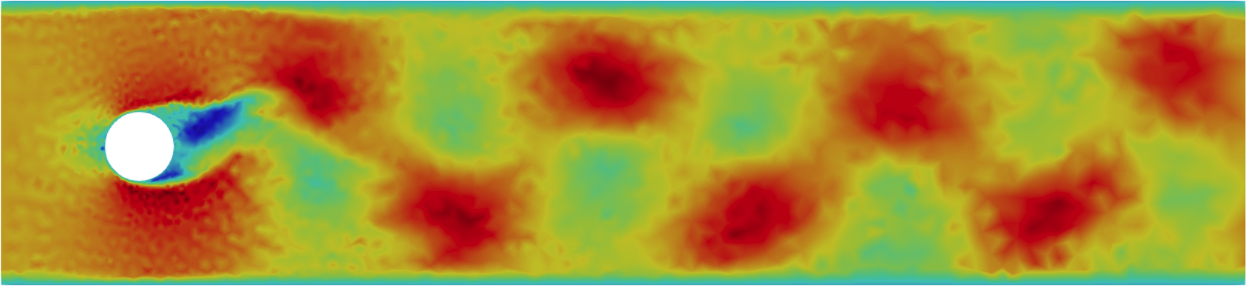}
\end{minipage} &
\begin{minipage}{0.45\linewidth}
\includegraphics[width=\linewidth]{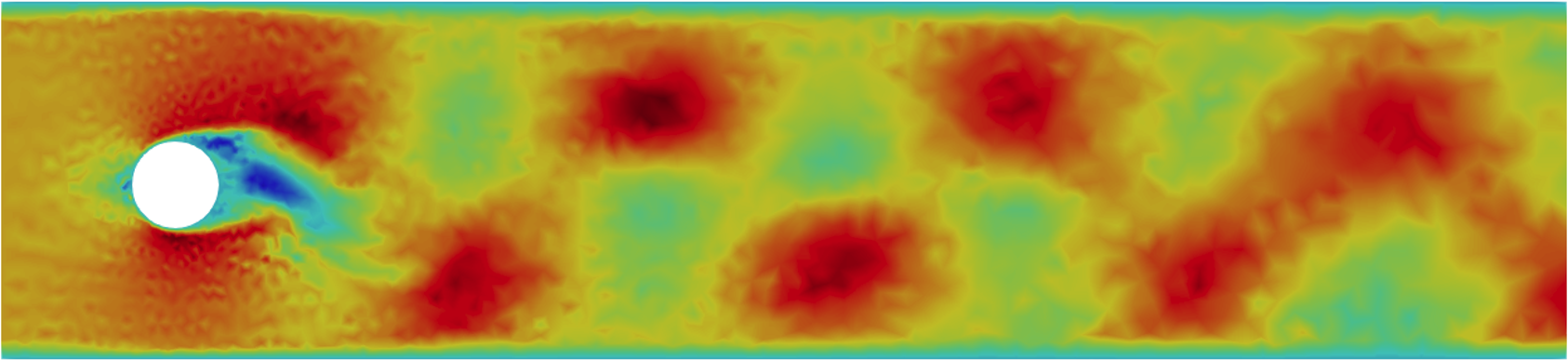}
\end{minipage}\\
(a) {\small Full model, $t = 14.6s$} & (b) {\small Full model, $t = 14.6s$} \\

\begin{minipage}{0.45\linewidth}
\includegraphics[width=\linewidth]{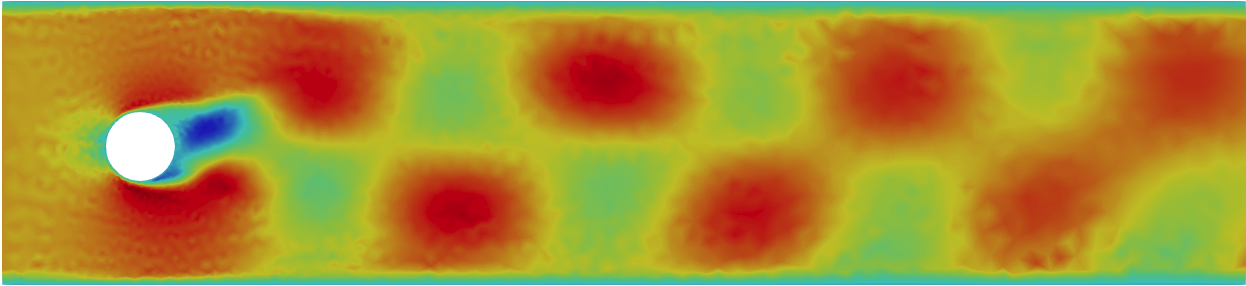}
\end{minipage} &
\begin{minipage}{0.45\linewidth}
\includegraphics[width=\linewidth]{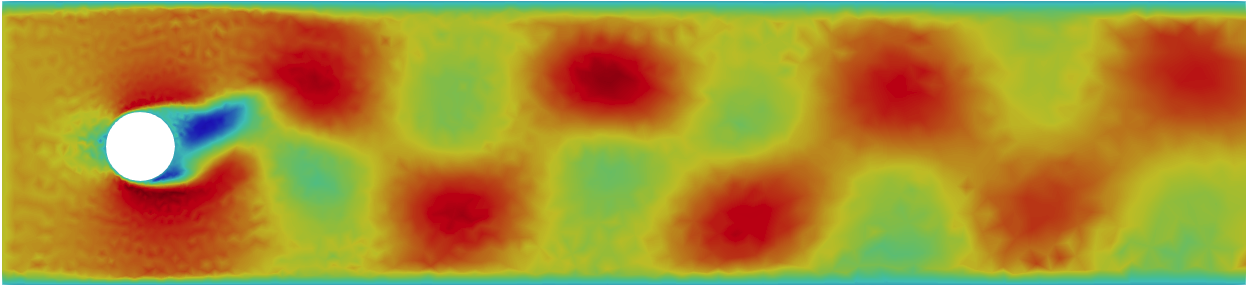}
\end{minipage} \\
(c) {\small POD+LSTM, $r = 2, t = 14.6$s} & (d) {\small POD+LSTM, $r = 6, t = 14.6s$} \\

\begin{minipage}{0.45\linewidth}
\includegraphics[width=\linewidth]{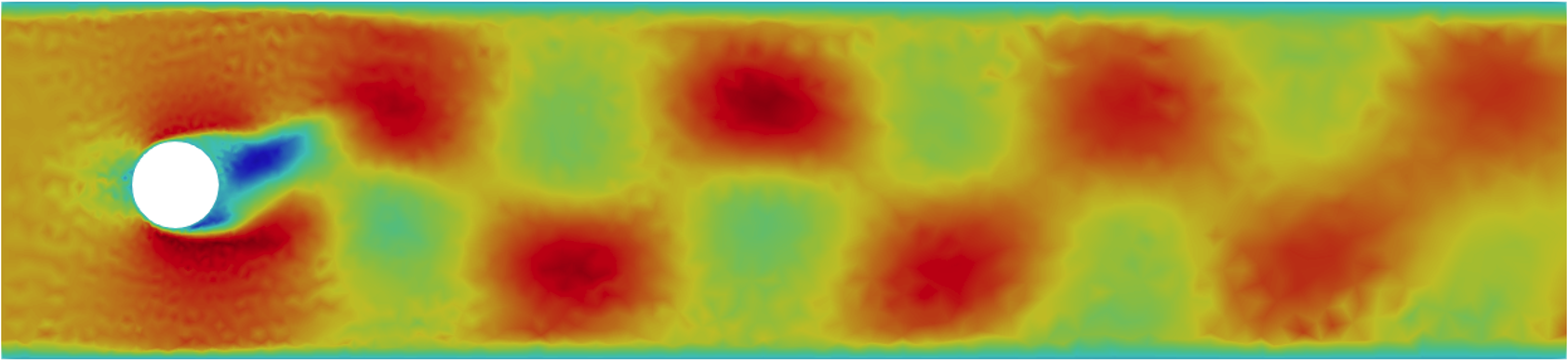}
\end{minipage} &
\begin{minipage}{0.45\linewidth}
\includegraphics[width=\linewidth]{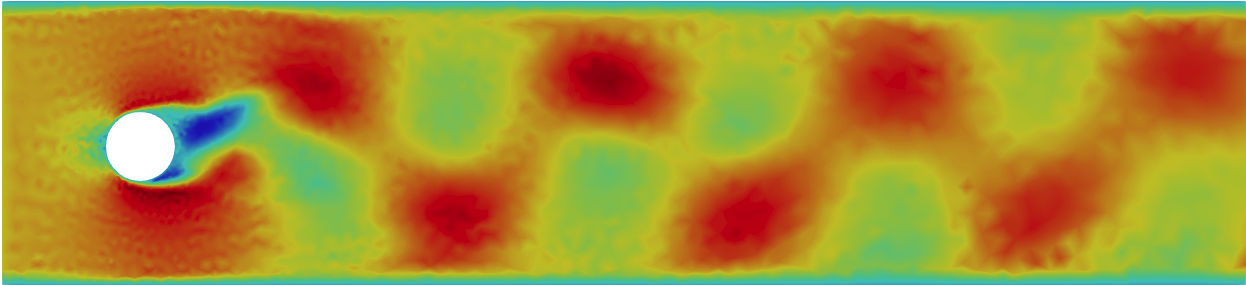}
\end{minipage} \\
(e) {\small PMD, $r = 2, t = 14.6$s} & (f) {\small PMD, $r = 6, t = 14.6s$} \\

\begin{minipage}{0.45\linewidth} 
\includegraphics[width = \linewidth,angle=0,clip=true]{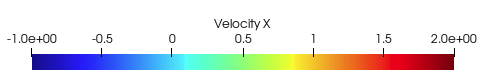}
\end{minipage}
&
\begin{minipage}{0.45\linewidth} 
\includegraphics[width = \linewidth,angle=0,clip=true]{scale1.png}
\end{minipage}\\

\end{tabular}
\caption{Flow past a cylinder test case (Re=5000): velocity solutions at $t = 14.6s$ obtained from the high fidelity full model, POD+LSTM and PMD using 2 and 6 basis functions respectively.}
\label{fg:flowconstruct}
\end{figure}

The computational domain is \(1.8 \times 0.41\) units, and a cylinder of radius 0.05 centered at \((0.2, 0.2)\). Fluid enters through the left boundary with a uniform inflow velocity of 1, navigates around the cylinder, and exits through the right boundary. No-slip and zero outflow conditions are applied to the top and bottom boundaries, while Dirichlet boundary conditions are enforced on the cylinder's surface.

The simulation time period is [0-20] seconds, and the computational mesh has 3802 nodes, and 200 snapshots captured at regular intervals of \(\Delta t = 0.1\). 130 snapshots (13 seconds) are used for constructing the ROM, and the remaining 70 snapshots are used to test the predictive capabilities of the PMD based ROM.

\begin{figure}[tbhp]
\centering
\begin{tabular}{cc}
\begin{minipage}{0.45\linewidth}
\includegraphics[width=\linewidth]{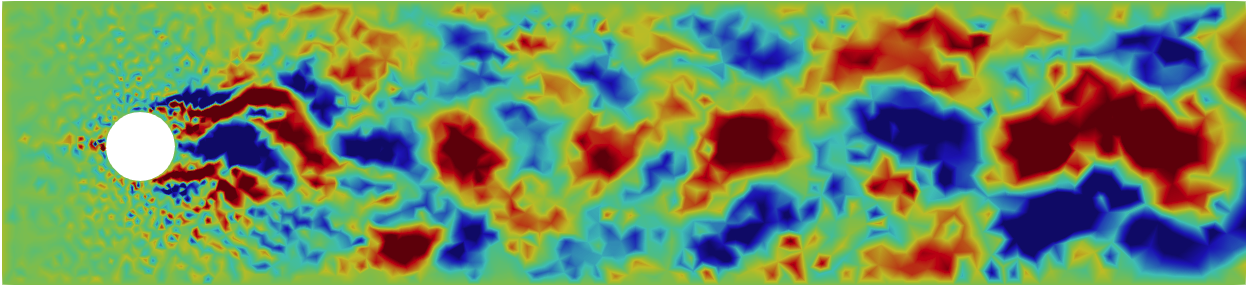}
\end{minipage} &
\begin{minipage}{0.45\linewidth}
\includegraphics[width=\linewidth]{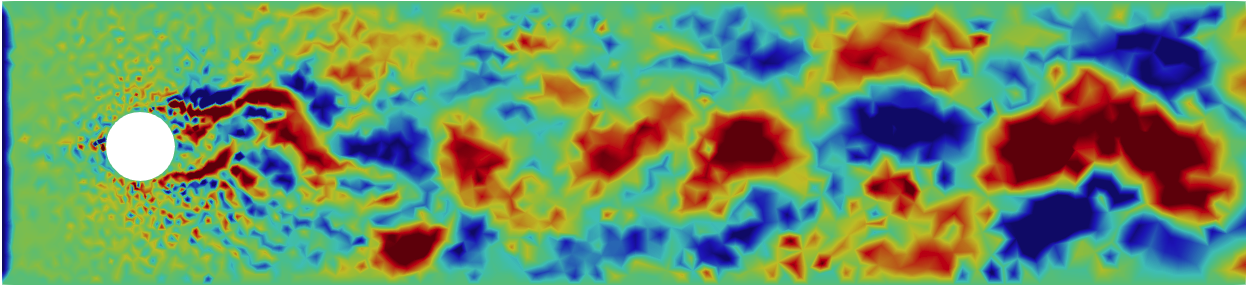}
\end{minipage} \\
(a) {\small POD+LSTM, $r = 2, t = 14.6$s} & (b) {\small POD+LSTM, $r = 6, t = 14.6$s} \\

\begin{minipage}{0.45\linewidth}
\includegraphics[width=\linewidth]{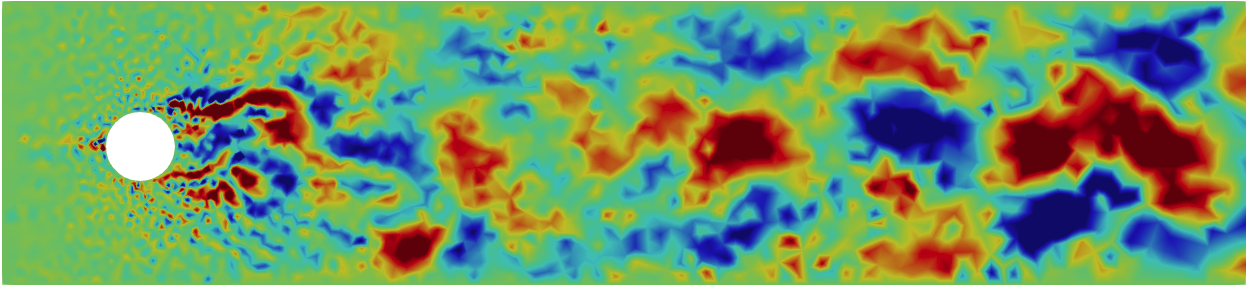}
\end{minipage} &
\begin{minipage}{0.45\linewidth}
\includegraphics[width=\linewidth]{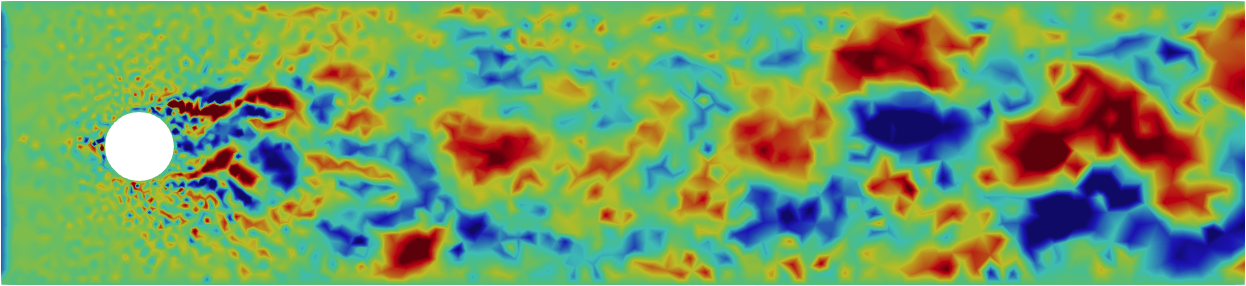}
\end{minipage} \\
(c) {\small PMD, $r = 2, t = 14.6$s} & (d) {\small PMD, $r = 6, t = 14.6$s} \\

\begin{minipage}{0.45\linewidth} 
\includegraphics[width = \linewidth,angle=0,clip=true]{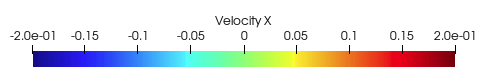}
\end{minipage}
&
\begin{minipage}{0.45\linewidth} 
\includegraphics[width = \linewidth,angle=0,clip=true]{scale2.png}
\end{minipage}\\

\end{tabular}
\caption{Flow past a cylinder test case (Re=5000): errors of POD+LSTM and PMD at $t = 14.6s$ using 2 and 6 basis functions respectively.}
\label{fg:flowerror1}
\end{figure}

\begin{figure}[tbhp]
\centering
\begin{tabular}{cc}
\begin{minipage}{0.45\linewidth}
\includegraphics[width=\linewidth]{flowfullmodel14.6s.png}
\end{minipage} &
\begin{minipage}{0.45\linewidth}
\includegraphics[width=\linewidth]{flowfullmodel14.9s.png}
\end{minipage} \\
(a) {\small Full model, $t = 14.6s$} & (b) {\small Full model, $t = 14.9s$} \\

\begin{minipage}{0.45\linewidth}
\includegraphics[width=\linewidth]{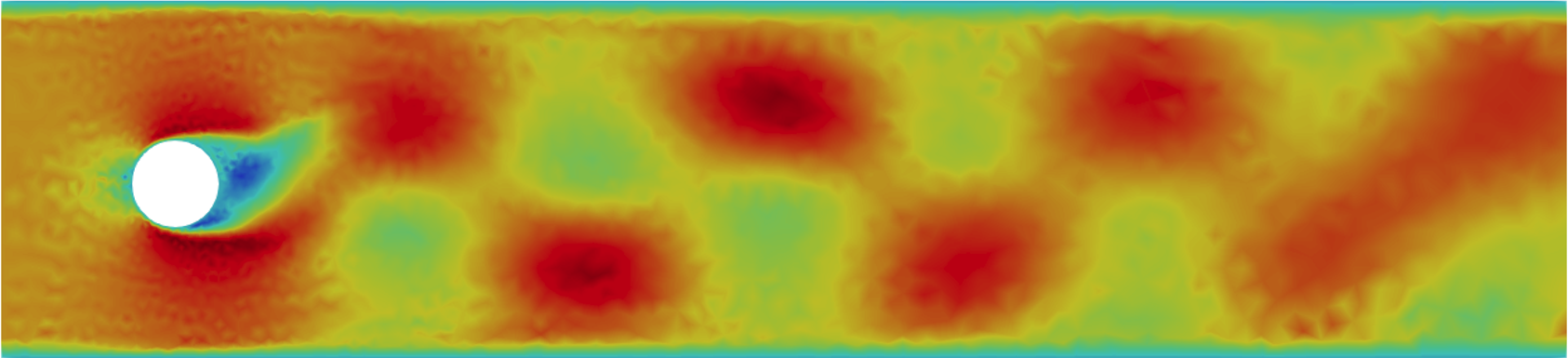}
\end{minipage} &
\begin{minipage}{0.45\linewidth}
\includegraphics[width=\linewidth]{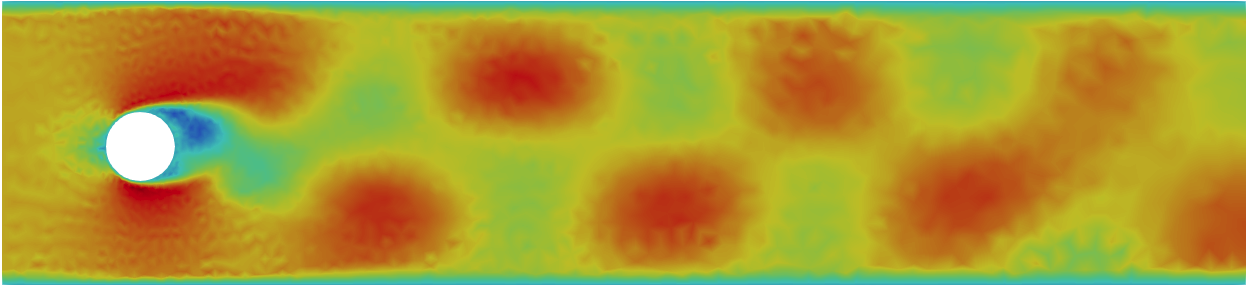}
\end{minipage} \\
(c) {\small POD+LSTM, $r = 6, t = 14.6s$} & (d) {\small POD+LSTM, $r = 6, t = 14.9s$} \\

\begin{minipage}{0.45\linewidth}
\includegraphics[width=\linewidth]{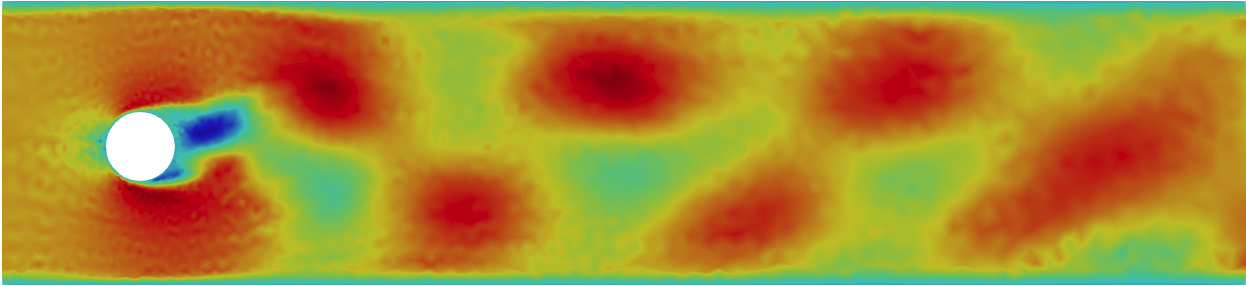}
\end{minipage} &
\begin{minipage}{0.45\linewidth}
\includegraphics[width=\linewidth]{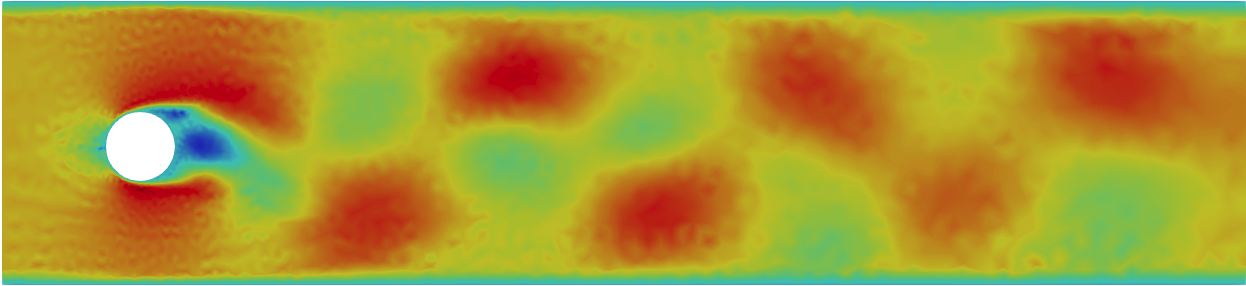}
\end{minipage} \\
(e) {\small POD+LSTM, $r = 12, t = 14.6s$} & (f) {\small POD+LSTM, $r = 12, t = 14.9s$} \\

\begin{minipage}{0.45\linewidth}
\includegraphics[width=\linewidth]{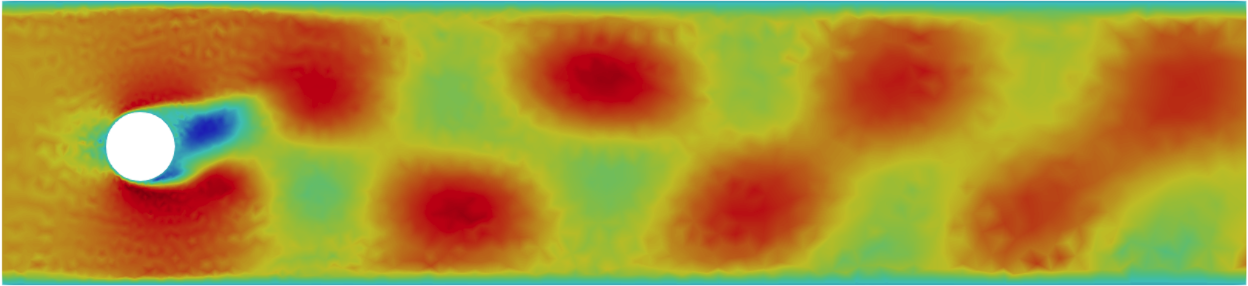}
\end{minipage} &
\begin{minipage}{0.45\linewidth}
\includegraphics[width=\linewidth]{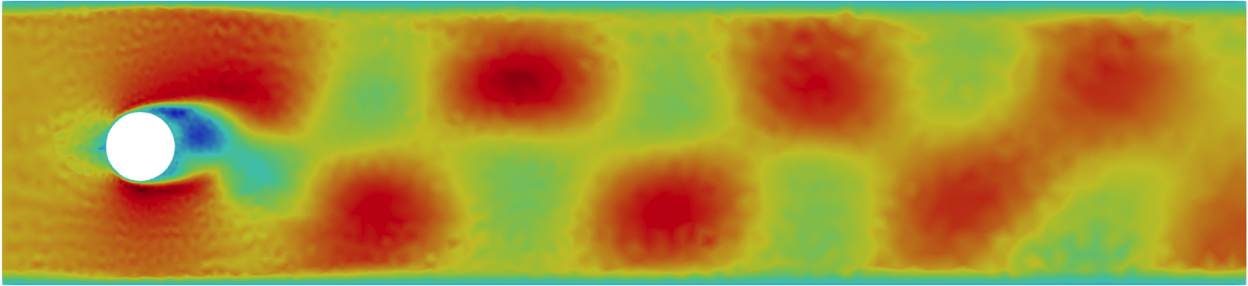}
\end{minipage} \\
(g) {\small PMD, $r = 6, t = 14.6s$} & (h) {\small PMD, $r = 6, t = 14.9s$} \\

\begin{minipage}{0.45\linewidth}
\includegraphics[width=\linewidth]{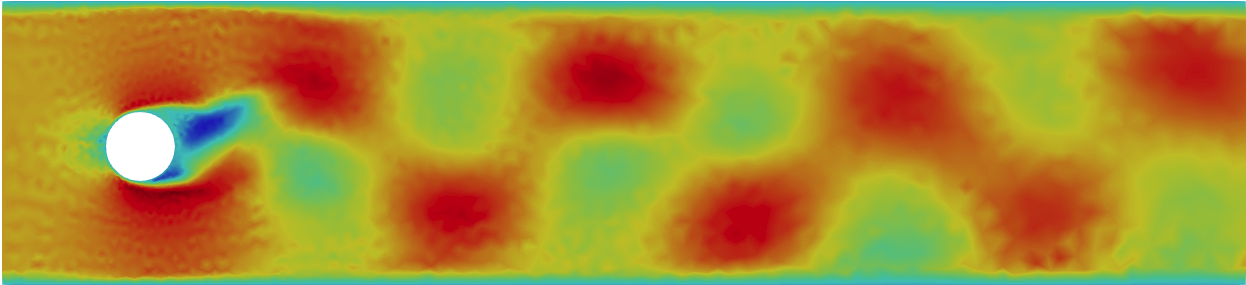}
\end{minipage} &
\begin{minipage}{0.45\linewidth}
\includegraphics[width=\linewidth]{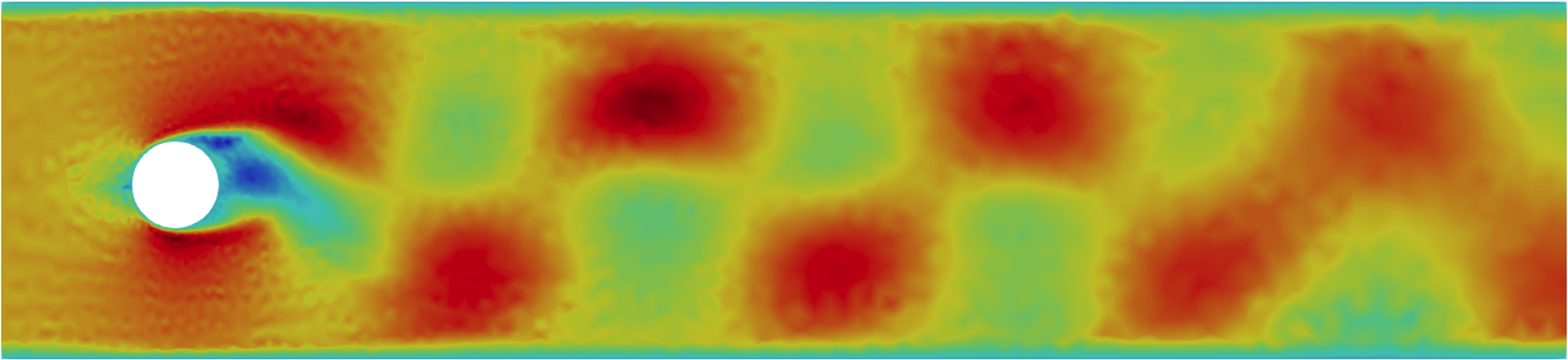}
\end{minipage} \\
(i) {\small PMD, $r = 12, t = 14.6s$} & (j) {\small PMD, $r = 12, t = 14.9s$} \\

\begin{minipage}{0.45\linewidth} 
\includegraphics[width = \linewidth,angle=0,clip=true]{scale1.png}
\end{minipage}
&
\begin{minipage}{0.45\linewidth} 
\includegraphics[width = \linewidth,angle=0,clip=true]{scale1.png}
\end{minipage}\\

\end{tabular}
\caption{Flow past a cylinder test case: velocity solutions obtained from the high fidelity full model, POD+LSTM and PMD at $t = 14.6s$ and $t=14.9s$ using 6 and 12 basis functions respectively.}
\label{fg:flowpredict}
\end{figure}

\begin{figure}[tbhp]
\centering
\begin{tabular}{cc}
\begin{minipage}{0.45\linewidth}
\includegraphics[width=\linewidth]{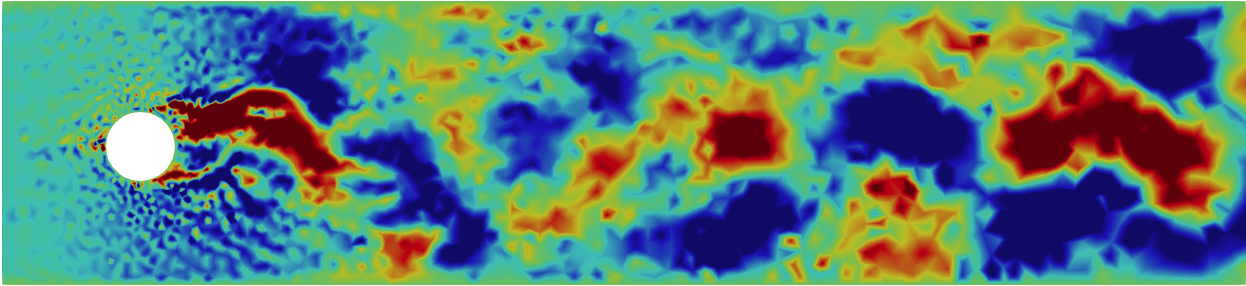}
\end{minipage} &
\begin{minipage}{0.45\linewidth}
\includegraphics[width=\linewidth]{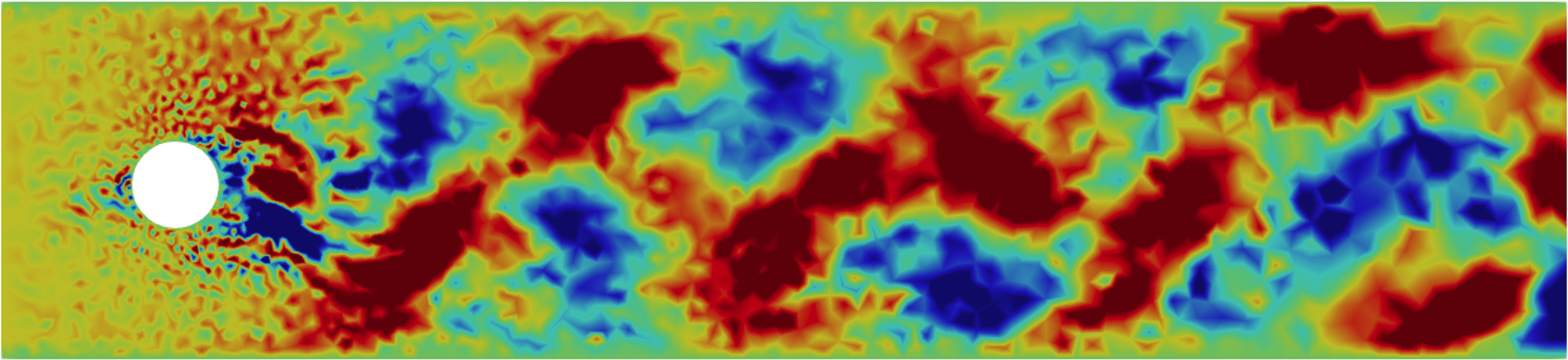}
\end{minipage} \\
(a) {\small POD+LSTM, $r = 6, t = 14.6s$} & (b) {\small POD+LSTM, $r = 6, t = 14.9s$} \\

\begin{minipage}{0.45\linewidth}
\includegraphics[width=\linewidth]{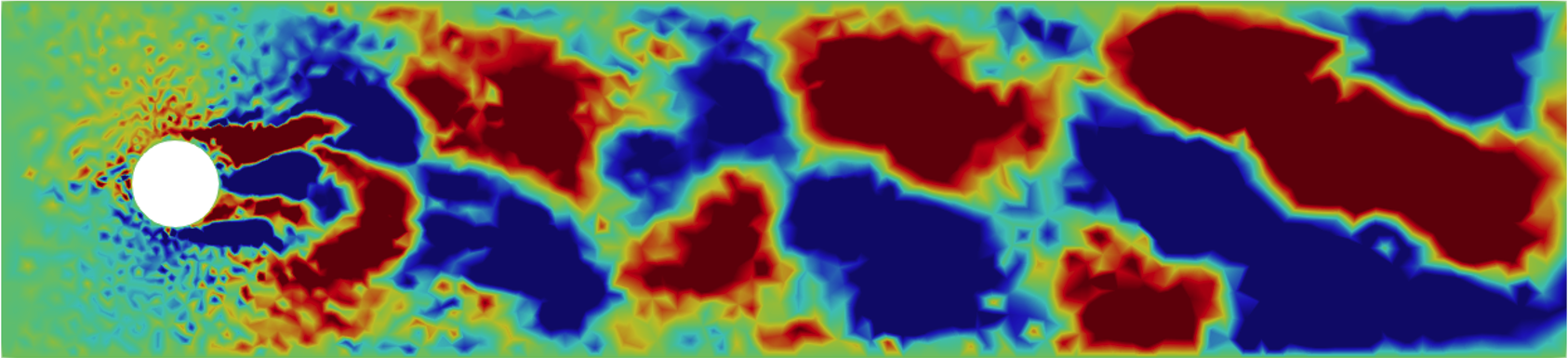}
\end{minipage} &
\begin{minipage}{0.45\linewidth}
\includegraphics[width=\linewidth]{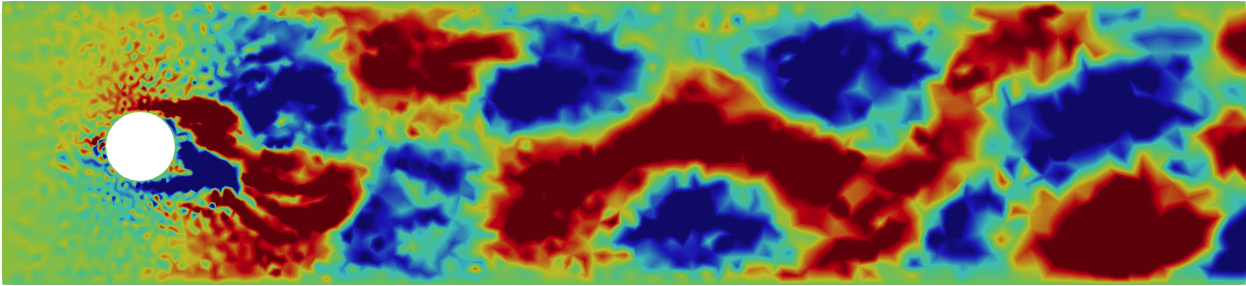}
\end{minipage} \\
(c) {\small POD+LSTM, $r = 12, t = 14.6s$} & (d) {\small POD+LSTM, $r = 12, t = 14.9s$} \\

\begin{minipage}{0.45\linewidth}
\includegraphics[width=\linewidth]{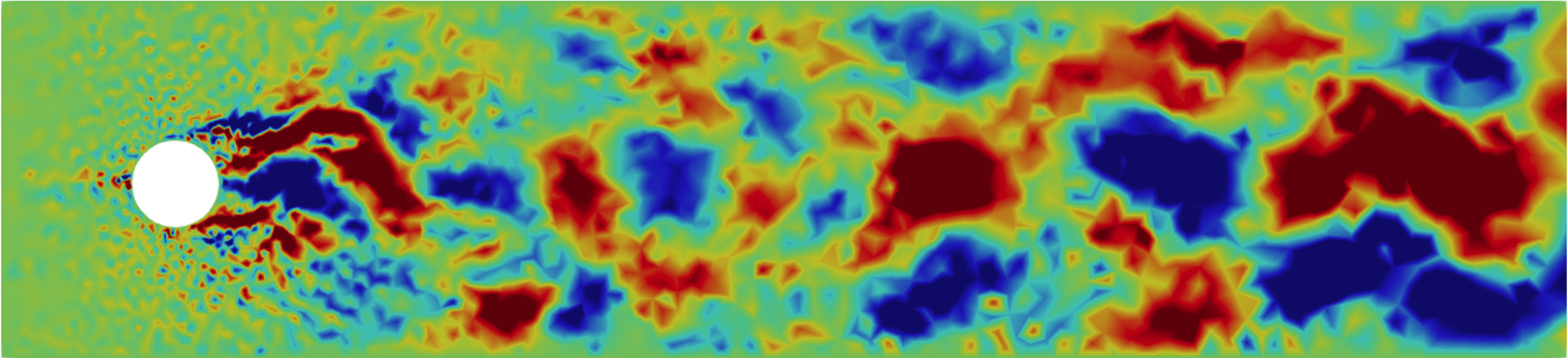}
\end{minipage} &
\begin{minipage}{0.45\linewidth}
\includegraphics[width=\linewidth]{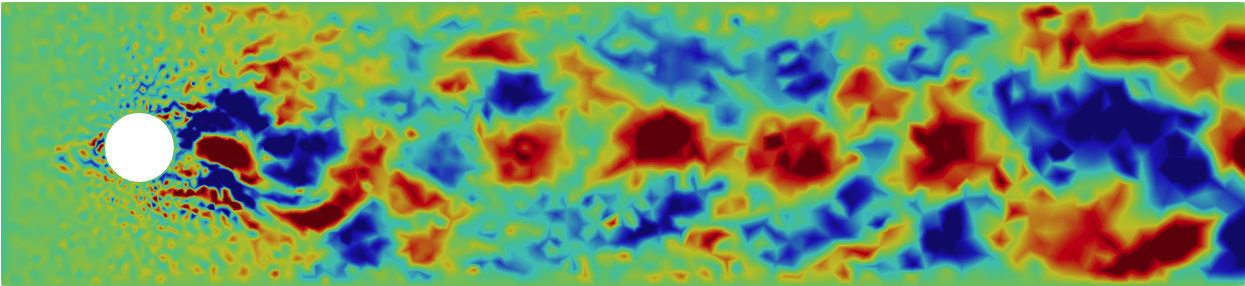}
\end{minipage} \\
(e) {\small PMD, $r = 6, t = 14.6s$} & (f) {\small PMD, $r = 6, t = 14.9s$} \\

\begin{minipage}{0.45\linewidth}
\includegraphics[width=\linewidth]{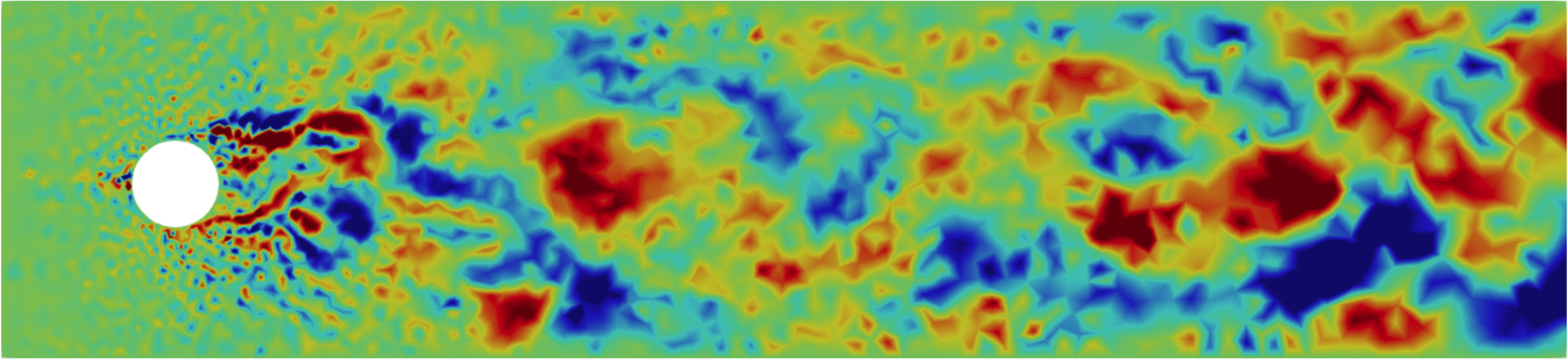}
\end{minipage} &
\begin{minipage}{0.45\linewidth}
\includegraphics[width=\linewidth]{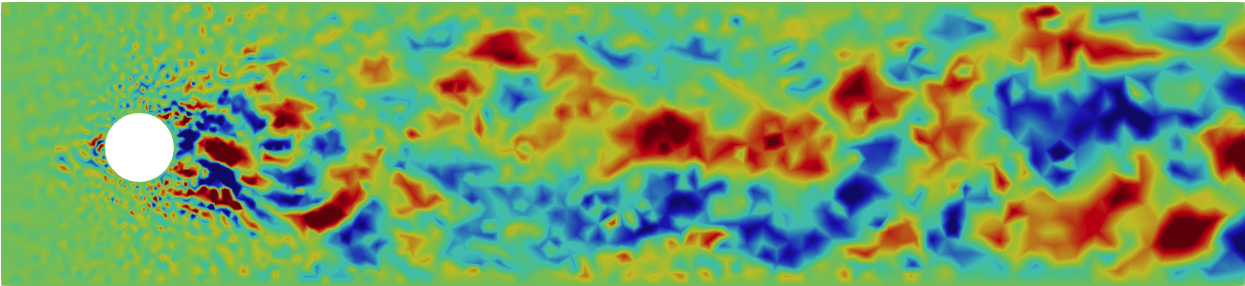}
\end{minipage} \\
(g) {\small PMD, $r = 12, t = 14.6s$} & (h) {\small PMD, $r = 12, t = 14.9s$} \\

\begin{minipage}{0.45\linewidth} 
\includegraphics[width = \linewidth,angle=0,clip=true]{scale2.png}
\end{minipage}
&
\begin{minipage}{0.45\linewidth} 
\includegraphics[width = \linewidth,angle=0,clip=true]{scale2.png}
\end{minipage}\\

\end{tabular}
\caption{Flow past a cylinder test case: errors of POD+LSTM and PMD at $t = 14.6s, 14.9s$ using 6 and 12 basis functions.}
\label{fg:flowerror2}
\end{figure}

\ref{fg:flowconstruct} shows the velocity results of the flow past a cylinder (Re = 5000) test case obtained from the high fidelity full model, PMD and ROM based on POD and LSTM using 2 and 6 basis functions respectively. The errors are also compared in \ref{fg:flowerror1}. As shown in the figures, POD struggles to capture nonlinear flow features using 2 basis functions, particularly failing to capture wake vortices and fine-scale dynamics. Increasing the number of basis functions to 6 improves the accuracy (see \ref{fg:flowconstruct} (d)), but some small-scale details still cannot be captured. \ref{fg:flowerror1} shows that the errors of PMD using 6 basis functions are smaller than those of PMD using 2 basis functions and POD using 6 basis functions. 
These results highlight that PMD performs better than the POD based ROM using the same number of basis functions.

In order to test the predictive capabilities of PMD outside the data range of constructing the ROM, results from two unseen time levels ($t=14.6$s and $t=14.9$s) are given. 
\ref{fg:flowpredict} and \ref{fg:flowerror2} compare the predictive performance of POD+LSTM and PMD for the flow past a cylinder problem (Re = 5000) using 6 and 12 basis functions. As shown in \ref{fg:flowpredict}, PMD consistently outperforms POD+LSTM for predicting next time levels using the same number of basis functions. The prediction errors of both methods are analyzed in \ref{fg:flowerror2}. The results indicate that increasing the number of basis functions improves the accuracy of both methods. However, PMD achieves better accuracy than POD+LSTM using the same number of basis functions or modes. 


\section{Conclusions}  \label{sec:summary}
A new non-linear, non-intrusive reduced-order model, probabilistic manifold decomposition (PMD) method, is presented in this work. The PMD consists of two main components: linear dimensionality reduction using Singular Value Decomposition (SVD) and the construction of a low-dimensional probabilistic manifold via a Markov process and geodesic distance learning. The SVD part is used to capture the main features of the system, while the probabilistic manifold part captures the remaining complex non-linear dynamics of the system by embedding the residual into a low-dimensional manifold.

Numerical experiments were conducted to validate the effectiveness of the PMD method using Fluidity \cite{fluidity_manual}, which is an open source finite element CFD model. Two test cases were considered: flow past a cylinder case with a Reynolds number $Re=5000$, and the lock exchange problem. The results show that PMD not only reduces the computational cost drastically but also preserves the accuracy of the system. 
Unlike recent deep learning methods such as POD combined with LSTM, PMD does not take a long time to train, and it is not like a black box. The PMD has explicit equations to explain the system. 

Another advantage of PMD is that it efficiently reduces dimensionality while capturing the high non-linearity in the systems, making it suitable for  problems with highly non-linearity. The PMD method will be extended to more complicated fluid dynamics applications in the future, such as multi-phase flow, large-scale urban flows, and ocean modeling. The predictive capability of PMD is closely linked to the amount of training data and the quality of the probabilistic manifold construction.

\section*{Acknowledgements}
\noindent 
The authors acknowledge the support of the Fundamental Research Funds for the Central Universities, the Top Discipline Plan of Shanghai Universities-Class I and Shanghai Municipal Science and Technology Major Project (No. 2021SHZDZX0100), National Key $R\&D$ Program of China(NO.2022YFE0208000). This work is supported in part by grants from the Shanghai Engineering Research Center for Blockchain Applications And Services (No. 19DZ2255100) and the Shanghai Institute of Intelligent Science and Technology, Tongji University.	

\bibliographystyle{unsrt} 
\bibliography{references}

\end{document}